\documentclass[pdflatex,sn-mathphys-num]{sn-jnl}

\usepackage{graphicx}%
\usepackage{multirow}%
\usepackage{amsmath,amssymb,amsfonts}%
\usepackage{amsthm}%
\usepackage{mathrsfs}%
\usepackage[title]{appendix}%
\usepackage{xcolor}%
\usepackage{textcomp}%
\usepackage{manyfoot}%
\usepackage{booktabs}%
\usepackage{algorithm}%
\usepackage{algorithmicx}%
\usepackage{algpseudocode}%
\usepackage{listings}%
\usepackage{bm} 

\theoremstyle{thmstyleone}%
\newtheorem{theorem}{Theorem}[section]
\newtheorem{lemma}[theorem]{Lemma}%

\theoremstyle{thmstyletwo}%
\newtheorem{remark}{Remark}%

\theoremstyle{thmstylethree}%

\begin{document}

\title[]{Non-convex optimal control governed by nonlinear elliptic partial differential equation: Analysis, functional-based algorithm and control landscape}


\author[1]{\fnm{Ning} \sur{Du}}\email{duning@sdu.edu.cn}

\author[1]{\fnm{Yanlin} \sur{Liu}}\email{yanlinliu@mail.sdu.edu.cn}

\author[2]{\fnm{Lei} \sur{Zhang}}\email{zhangl@math.pku.edu.cn}

\author[1]{\fnm{Xiangcheng} \sur{Zheng}}\email{xzheng@sdu.edu.cn}

\affil[1]{\orgdiv{School of Mathematics}, \orgname{Shandong University}, \orgaddress{\city{Jinan}, \postcode{250100}, \country{China}}}

\affil*[2]{\orgdiv{School of Mathematical Sciences, Beijing International Center for Mathematical Research,
  Center for Machine Learning Research,
  Center for Quantitative Biology}, \orgname{Peking University}, \orgaddress{\city{Beijing}, \postcode{100871}, \country{China}}}


\abstract{
Non-convex optimal control arises from various applications but may contain multiple stationary points.
Classical solvers usually perform a local search and therefore rely on good initial guesses to reach appropriate local optimal controls. In this work we introduce a novel solution strategy for the non-convex optimal control of an elliptic equation, the main idea of which is to construct the control landscape based on the functional high-index saddle dynamics (FHiSD) method.
This method reduces the dependence on prescribed initial guesses by using the information of high-index saddle points,  and depicts the macroscopic configuration of the space of the control variable. Then various minima could be systematically computed along transition pathways and the control strategy can then be selected among them. 
We prove the existence of the optimal control, and then propose and analyze the FHiSD in locating saddle points. Subsequently, the applicability of FHiSD in non-convex optimal control is rigorously justified.
Numerical results not only indicate the effectiveness of the proposed method, but reveal unintuitive phenomena (e.g. the non-monotonicity of the values of the cost functional with respect to the Morse indices)  that support the necessity of computing multiple solutions of high indices.}

\keywords{non-convex, optimal control, saddle dynamics, control landscape, multiple solutions}


\pacs[MSC Classification]{49M41, 65K10}

\maketitle

\section{Introduction}

\subsection{Motivations}
Non-convex optimal control is widespread in scientific and engineering problems 
\cite{friedemannFindingGlobalSolutions2023,liuPosterioriErrorEstimates2003}, and extensive numerical methods have been developed, see e.g. the gradient and Newton-type methods \cite{NumericalOptimization2006a,troltzschOptimalControlPartial2010}, 
the sequential quadratic programming \cite{OptimizationPDEConstraints2009,troltzschOptimalControlPartial2010}, and the penalized gradient projection method \cite{chryssoverghiDescentPenaltyMethodsRelaxed2008,coletsosOptimalControlNonlinear2014}. Nevertheless, most methods perform a local study without considering the macroscopic configuration. Specifically, they usually focus on investigations of a saddle point or local minimum of the non-convex optimal control, and thus rely on good initial conditions to reach appropriate control strategies.
To address these issues, this work intends to propose a novel solution strategy based on systematic saddle point algorithms, which investigates the non-convex optimal control from a comprehensive perspective under reduced requirements on initial guesses.  

For clarity, we first introduce several concepts.
Let $X$ be a Hilbert space, and let $J:X\to\mathbb{R}$ be a three times Fr\'echet differentiable energy functional. 
A stationary point $x^*\in X$ of $J$ satisfies $\nabla J(x^*)=0$. If $x^*$ is not a local minimum, then we call $x^*$ a saddle point. The Morse index of a nondegenerate saddle point is the maximal dimension of a subspace of $X$ on which the Hessian is negative definite \cite{milnor2016morse}, and we call a saddle point with Morse index $k$ an index-$k$ saddle point (a minimum is indeed an index-$0$ saddle point).

There exist several efficient methods for locating saddle points \cite{farrellComputingMultipleSolutions2026,MR3522976,22M1478598,LiuXie,quappLocatingSaddlePoints2014,yinHighIndexOptimizationBasedShrinking2019a}, one of which is the  high-index saddle dynamics (HiSD) method. 
With the help of HiSD, the solution landscape method, which aims to construct a pathway map containing saddle points and local minima of a system and the connections among them,  has been developed in \cite{yinSearchingSolutionLandscape2021b} and successfully applied in several fields, cf. \cite{wangModellingComputationLiquid2021,zhang2026cellsystems,wu2026jamming}. 
A typical procedure to construct the solution landscape starts from a high-index saddle point (i.e. the parent state) and then locates lower-index saddle points by the HiSD \cite{yinSearchingSolutionLandscape2021b}. 
Since the initial search directions of HiSD for a certain saddle point are provided from eigenvectors of the Hessian at a higher-index saddle point, the solution landscape method circumvents the difficulty of selecting proper initial values and makes the downward search for lower-index saddle points efficient. 

The above observations naturally motivate the application of the HiSD-based solution landscape method to the non-convex optimal control problems, while substantial difficulties are encountered immediately. 
Specifically, the traditional HiSD is mainly formulated in Euclidean spaces and relies on explicit gradient and Hessian expressions of the model. In contrast, optimal control of partial differential equations is proposed in Sobolev spaces and it may be difficult to give explicit gradient and Hessian expressions of the corresponding reduced cost functionals due to, e.g. the state equation constraint and the $H^1$ regularization.
These difficulties motivate us to develop novel HiSD methods that are naturally compatible with optimal control problems.

\subsection{Non-convex optimal control model}

Elliptic optimal control problems and their numerical discretizations have been extensively studied; see, e.g., \cite{18M1227147,casasOptimalControlPartial2017,OptimizationPDEConstraints2009,troltzschOptimalControlPartial2010}. 
In this work, we demonstrate the novel solution strategy of non-convex optimal control by the following problem
\vspace{-2pt}
\begin{equation}\label{equ:nc-ocp}
    \min\limits_{u\in H^1_0(\Omega)} J(y, u) := \frac{1}{2} \|y - y_d\|_{L^2(\Omega)}^2 
+ \frac{\lambda}{2} \|u\|_{H^1(\Omega)}^2,
\end{equation}
 subject to a general elliptic equation with nonlinear mappings on both $y$ and $u$ (cf.  \cite{betzOptimalControlNonsmooth2025,casasOptimalControlProblems2020} for similar nonlinear state equations in optimal control)
\vspace{-2pt}
\begin{equation}\label{equ:se-nc-ocp}
\mathcal{A} y + d(x,y) :=-\sum_{i,j=1}^{n} \partial_{x_i} \big(  a_{ij}(x) \partial_{x_j} y  \big) + c(x)y+d(x,y)= g(x,u)  \text{ in } \Omega;~
y = 0  \text{ on } \Gamma.
\end{equation}
Here $L^q(\Omega)$ with $1\le q\le\infty$ and $H^m(\Omega):=W^{m,2}(\Omega)$ with integer $m\ge0$ denote the standard Lebesgue and Sobolev spaces, respectively, and $H_0^1(\Omega)$ denotes the closure of $C_c^\infty(\Omega)$ in $H^1(\Omega)$. The $\Omega \subset \mathbb{R}^n$ ($1\leq n\leq 3$) is a bounded domain with a ${C}^{1,1}$ boundary $\Gamma$, and $y$ and $u$ denote the state and control variables, respectively.
The $y_d(x)$ denotes the desired state and $\lambda > 0$ is the regularization parameter. In the state equation, $\mathcal{A}$ is a symmetric and coercive second-order elliptic operator
where $ a_{ij}(x)\in C^1(\overline{\Omega})$, the space of continuously differentiable functions on $\overline{\Omega}$, $c(x) \in {L}^\infty (\Omega)$,  and
$d(x,y)$ and $g(x,u)$ are given functions. 


\subsection{Contributions of the work}
The main contributions of the work are enumerated as follows:

\begin{itemize}

\item For the optimal control problem \eqref{equ:nc-ocp}--\eqref{equ:se-nc-ocp}, we establish the well-posedness of the state equation and the existence of an optimal control. Compared with existing studies for elliptic optimal control problems, additional analysis is required to accommodate the nonlinear mapping of the control variable. In particular, we establish the mapping properties and local Lipschitz continuity of the associated Nemytskii and control-to-state operators under the nonlinear control mapping, which are used to justify the compactness and convergence arguments in the proof of the existence.

\item We develop a functional high-index saddle dynamics (FHiSD) in the $H^{-1}(\Omega)$ sense for locating saddle points of functionals defined on $H_0^1(\Omega)$. Concerning the impacts of the unbounded differential operators and the off-diagonal terms in the linearized system, we perform the stability analysis of FHiSD by means of solution operators. The analysis differs from that of the conventional HiSD in \cite{yinHighIndexOptimizationBasedShrinking2019a} , which is finite-dimensional and thus does not suffer the aforementioned issues. We further justify the applicability of FHiSD to the non-convex optimal control by analyzing the spectral structure of its Hessian and the boundedness of high-order derivative terms.

\item  We discretize the FHiSD by the finite element method to construct control landscapes. Compared with traditional local solvers, the proposed strategy systematically locates multiple stationary points with reduced dependence on initial guesses, while the numerical results reveal a non-monotonic relation between reduced cost values and saddle indices, highlighting the necessity of computing saddle-point solutions when seeking optimal controls.

\end{itemize}

The rest of the paper is organized as follows: 
In Section \ref{sec:optimal control}, we present a detailed analysis of the optimal control problem (\ref{equ:nc-ocp})--(\ref{equ:se-nc-ocp}).
In Section \ref{sec3}, we develop and analyze an FHiSD, and in Section \ref{sec:hisd} we justify its applicability for optimal control problems and discuss the construction of the control landscape.
In Section \ref{sec:experiments}, we perform several numerical experiments and comparisons. Some concluding remarks are addressed in the last section. 

\section{Analysis of optimal control}
\label{sec:optimal control}

This section presents the analysis of the non-convex optimal control problem (\ref{equ:nc-ocp})--(\ref{equ:se-nc-ocp}).
We assume that the semilinear term 
$d(x,y)$ satisfies standard local Lipschitz continuity and monotonicity, see e.g. the book of Tr\"oltzsch \cite[Assumption 4.14]{troltzschOptimalControlPartial2010}. We also
impose local or global Lipschitz continuity (depending on the dimension $n$ of the problem) and growth conditions on 
$g(x,u)$ to ensure the well-posedness of the state equation. 

\vspace{0.05in}

\noindent\underline{For $n=1$, we impose Assumptions A and B as follows:}

\vspace{0.05in}

\noindent\textbf{Assumption A}\\
(i) ~$d(\cdot,\cdot):\Omega \times \mathbb{R} \to \mathbb{R}$ is a Carathéodory function.\\
(ii) $d(x,y)$ is three times continuously differentiable with respect to $y$. There exists a $K>0$ such that  
\begin{align}
\lvert  D_y^j d(x,0)\rvert \leq K, \text{ for }  0 \leq j \leq 2, \text{ and a.e. } x \in \Omega, \label{equ:bc-d} 
\end{align}
and for any $E> 0$ there exists a positive constant $L_{E}$ such that
\begin{align}
\lvert D_y^j d(x,y_1)-D_y^j d(x,y_2) \rvert \leq L_{E} \left\lvert y_1-y_2 \right\rvert ,\text{ for } 0 \leq j \leq 2 \text{ and a.e. } x \in \Omega,\label{equ:lip-d}
\end{align}
where $y_i \in \mathbb{R} $ satisfies $|y_i| \leq E$ for $ i=1,2$.\\
(iii) $D_y d(x,y) \ge 0 \text{ for any } y \in \mathbb{R}  \text{ and a.e. } x \in \Omega.$

\vspace{0.05in}

\noindent\textbf{Assumption B}\\
(i) ~$g(\cdot,\cdot):\Omega \times \mathbb{R} \to \mathbb{R}$ is a Carathéodory function.\\
(ii) $g(x,u)$ is three times continuously differentiable with respect to $u$, and for any $M > 0$, there exists a positive constant $L_{M}$ such that
\begin{equation}
\lvert D_u^j g(x,u_1)-D_u^j g(x,u_2) \rvert \leq L_{M} \left\lvert u_1-u_2 \right\rvert ,\text{ for } 0 \leq j \leq 2 \text{ and a.e. } x \in \Omega, \label{equ:lip-g1}
\end{equation}
where $u_i \in \mathbb{R} $ satisfies $|u_i| \leq M$ for $ i=1,2$.\\
(iii) There exist a constant $b>0$ and a function $a\in {L}^2(\Omega)$ such that $g(x,u)$ satisfies the growth condition for some $1 \leq s <\infty$
\begin{equation*} 
\left\lvert g(x,u) \right\rvert \leq a(x) + b\left\lvert u\right\rvert^{s}, \text{ for a.e. } x \in \Omega \text{ and any } u \in \mathbb{R}.
\end{equation*}

\noindent\underline{For $n=2,3$, we impose the Assumptions A and C:}

\vspace{0.05in}

\noindent\textbf{Assumption C}\\
(i) ~$g(\cdot,\cdot):\Omega \times \mathbb{R} \to \mathbb{R}$ is a Carathéodory function.\\
(ii) $g(x,u)$ is three times continuously differentiable with respect to $u$,  $g(\cdot,0)\in L^2(\Omega)$, and there exists an $L_0 > 0$ such that for any $u_1,u_2 \in \mathbb{R}$
\begin{equation*} 
\lvert D_u^j g(x,u_1)-D_u^j g(x,u_2) \rvert \leq L_{0} \left\lvert u_1-u_2 \right\rvert ,\text{ for } 0 \leq j \leq 2 \text{ and a.e. } x \in \Omega.
\end{equation*}

\begin{remark}
The dimension-dependent assumption on $g$ is caused by the regularity of the control. Since $u\in H_0^1(\Omega)$, the pointwise boundedness of $u$ follows from the Sobolev embedding $H_0^1(\Omega)\hookrightarrow L^\infty(\Omega)$ only for $n=1$. Hence the local Lipschitz condition on $g$ is sufficient in one dimension, whereas a global Lipschitz condition is imposed for $n=2,3$. This global condition can be relaxed if additional pointwise bounds on the control are available, for example under box constraints; see \cite{casasOptimalControlProblems2020}.
\end{remark}


In the rest of this section we analyze the optimal control problem (\ref{equ:nc-ocp})--(\ref{equ:se-nc-ocp}) for the one-dimensional case (i.e. $n=1$) under Assumptions A and B. The two- and three-dimensional problems can be analyzed in a similar and simpler manner under Assumptions A and C due to the global Lipschitz continuity of $g$, and we thus omit the details of the proofs. We will use $C$ to denote a generic positive constant that may assume different values at different occurrences, and we use $C(M)$ with some quantity $M$ to show the dependence of $C$ on $M$. Furthermore, we use $(\cdot,\cdot)$ and $\|\cdot\|$ to denote the  inner product and the norm of $L^2(\Omega)$, respectively, and denote $f_y:=D_y f$, $f_{yy}:=D_y^2 f$, and $f_{yyy}:=D_y^3 f$ for any function $f$ for simplicity.

\subsection{Well-posedness of state equation}
To facilitate the treatment of the nonlinear term $g(x,u)$, we introduce the associated Nemytskii operator $\mathcal{G}$ defined by $\mathcal{G}(u)(x) = g(x,u(x))$ for a.e. $x\in\Omega$. 
The following lemma describes some  properties of  $\mathcal{G}$.
\begin{lemma}\label{lm:mp-g}
Under the Assumption B for one-dimensional case (or Assumption C for two- and three-dimensional cases), $\mathcal{G} (u) \in {L} ^2(\Omega)$ for $u\in {H}_0^1(\Omega)$ and the estimate $\left\lVert \mathcal{G} (u) \right\rVert \leq C \big( 1+\left\lVert u \right\rVert_{{H}_0^1(\Omega)}^s \big)$ holds for some $ 1 \leq s <\infty $ given in Assumption B in the one-dimensional case or $s=1$ for two- and three-dimensional cases,
where $C$ is independent of $u$.

Furthermore, for any $M>0$,  the following estimate holds for $u_i \in {H}_0^1(\Omega) $ ($i=1,2 $) with $\left\lVert u_i \right\rVert_{{H}_0^1(\Omega)} \leq M$
\begin{equation} \label{equ:lip-g2}
    \left\lVert \mathcal{G}(u_1)-\mathcal{G}(u_2) \right\rVert \leq C(M) \left\lVert u_1-u_2 \right\rVert.
\end{equation}
\end{lemma}

\begin{proof}
We first consider the one-dimensional case. For every control $u\in {H}_0^1(\Omega)$ and a.e. $x \in \Omega$, we have
\begin{equation*}
\begin{aligned}
\left\lVert \mathcal{G} (u) \right\rVert^2   \leq \int_{\Omega} ( a(x) + b|u|^s )^2 \,dx  \leq C\big( 1 + \left\lVert u \right\rVert_{{L} ^{2s}(\Omega)}^{2s} \big)  \leq C\big( 1 + \left\lVert u \right\rVert_{{H}_0^1(\Omega)}^{2s} \big) , 
\end{aligned}
\end{equation*}
where the Sobolev embedding ${H}_0^1(\Omega) \hookrightarrow {L}^\infty(\Omega)  $ for $n=1$ is used. 
Given any $M>0$ and $u_i \in {H}_0^1(\Omega)$ satisfying $\left\lVert u_i \right\rVert_{{H}_0^1(\Omega)} \leq M$ for $i=1,2 $, their ${L} ^\infty(\Omega)$ norms are also bounded by some constant depending on $M$ such that we could apply the local Lipschitz continuity of $g(x,u)$ in (\ref{equ:lip-g1}) to get
\begin{equation*}
    \int_{\Omega} \left\lvert g(x,u_1) - g(x,u_2) \right\rvert^2 \,dx \leq  [C(M)]^2 \int_{\Omega} \left\lvert u_1 - u_2 \right\rvert^2 \,dx = [C(M)]^2 \left\lVert u_1 - u_2 \right\rVert^2.
\end{equation*}
The two- and three-dimensional cases can be proved in a similar and simpler manner without relying on $M$ based on the Assumption C, and we thus omit the details.
\end{proof}

Now we introduce the control-to-state operator $G:{H}_0^1(\Omega) \to Y:= {H}^2(\Omega) \cap {H}^1_0(\Omega)$. We first show that $G$ is well-defined.

\begin{theorem}\label{thm:wp-se}
Under the Assumptions A and B for one-dimensional case (or Assumptions A and C for two- and three-dimensional cases),
 the state equation (\ref{equ:se-nc-ocp}) with $u\in {H}_0^1(\Omega)$ has a unique  solution $y=G(u)\in Y$, which is continuous on $\overline{\Omega}$ and satisfies
$
\left\lVert y \right\rVert_{{H}^2(\Omega)} \leq  C\big( 1+\left\lVert u \right\rVert_{{H}_0^1(\Omega)}^s \big)$ for $ 1\leq s < \infty,
$
where $C$ is independent of $u$.
\end{theorem}

\begin{proof}
By Lemma \ref{lm:mp-g}, $u\in {H}_0^1(\Omega)$ implies $\mathcal{G}(u)\in {L}^2(\Omega)$ and then we could follow the same approach as \cite[Theorem 4]{casasOptimalControlPartial2017} to prove that the state equation (\ref{equ:se-nc-ocp}) has a unique solution $y \in Y$ that satisfies
$
\left\lVert y \right\rVert_{{H}^2(\Omega)} \leq C  \left\lVert 1+ \mathcal{G} (u)\right\rVert  \leq C\big( 1+\left\lVert u \right\rVert_{{H}_0^1(\Omega)}^s \big) .
$
The continuity of $y$ follows immediately from the Sobolev embedding theorem.
\end{proof}

Next, we show that $G$ is locally Lipschitz continuous with respect to $u$.
\begin{theorem} 
Under the Assumptions A and B for one-dimensional case (or Assumptions A and C for two- and three-dimensional cases),
 for any $M>0$  there exists a constant $C(M)>0$ such that the following relation holds for any $ u_i \in {H}_0^1(\Omega)$ with $\left\lVert u_i \right\rVert_{{H}_0^1(\Omega)} \leq M$ for $i =1,2$
\begin{equation}\label{equ:lip-G}
\left\lVert G(u_1)-G(u_2) \right\rVert_{{H}^2(\Omega)} \leq C(M) \left\lVert u_1-u_2\right\rVert.
\end{equation}
\end{theorem}

\begin{proof}
For any $M>0$, pick $u_1,u_2\in {H}_0^1(\Omega)$ such that $\left\lVert u_i \right\rVert_{{H}_0^1(\Omega)} \le M$ for $i =1,2$. By Theorem \ref{thm:wp-se} there exist continuous states $y_i=G(u_i)$ for $i=1,2$. 
We then subtract the corresponding state equations to get
\begin{equation*}
    \begin{aligned}
&\mathcal{A} (y_1-y_2) + d(x,y_1)-d(x,y_2) \\
&\quad=\mathcal{A} (y_1-y_2) + \int_0^1 d_y\big( x,y_2+\theta(y_1-y_2) \big)\, d\theta \cdot (y_1-y_2)\\
&\quad=:\mathcal{A} (y_1-y_2) + \phi(x)(y_1-y_2) = \mathcal{G}(u_1)-\mathcal{G}(u_2).
\end{aligned}
\end{equation*}
Thus, $y_1-y_2$ satisfies a linear elliptic equation with zero boundary conditions.
By Theorem \ref{thm:wp-se}, $\|u_i\|_{H_0^1(\Omega)}\le M$ and the Sobolev embedding $H^2(\Omega)\hookrightarrow L^\infty(\Omega)$, we have $\|y_i\|_{L^\infty(\Omega)}\le C(M)$ for $i=1,2$, which, together with the Assumption A, leads to $\|\phi\|_{L^\infty(\Omega)}\le C(M)$.
Moreover, $d_y\ge0$ implies $\phi(x)\ge0$.
Thus, the linear elliptic problem with respect to $y_1-y_2$ has a unique solution in ${H}^2(\Omega)$ with the regularity estimate \cite[Section 6.3, Theorem 4]{evans}
\begin{equation*}
\left\lVert y_1 -y_2 \right\rVert_{{H}^2(\Omega)}
\leq 
C(M) \left\lVert \mathcal{G}(u_1)-\mathcal{G}(u_2)\right\rVert 
\leq 
C(M) \left\lVert u_1-u_2\right\rVert 
\end{equation*} 
where we have used the local Lipschitz continuity of $\mathcal{G}$ in (\ref{equ:lip-g2}).
\end{proof}


\subsection{Existence of optimal control}
We prove the existence of an optimal control for the problem (\ref{equ:nc-ocp})--(\ref{equ:se-nc-ocp}) in the following theorem.

\begin{theorem} 
Under the Assumptions A and B for one-dimensional case (or Assumptions A and C for two- and three-dimensional cases),
there exists an optimal control $\overline{u} \in  {H}_0^1(\Omega)$ for  (\ref{equ:nc-ocp})--(\ref{equ:se-nc-ocp}) satisfying $J(\overline{y},\overline{u})\le J(y,u)$ for all $u\in {H}_0^1(\Omega)$, where $\overline{y} = G(\overline{u}) $ and $y = G(u) $.
\end{theorem}
\begin{proof}
Since $J$ is nonnegative and bounded below, there exists a minimizing sequence \cite[Theorem 2.16]{troltzschOptimalControlPartial2010}
$
    \{u_j\}_{j\in \mathbb{N}} \subset U:=\left\{u\in {H}_0^1(\Omega)| \left\lVert  u \right\rVert_{{H}_0^1(\Omega)} \leq B\right\} 
$
for some $B>0$ such that
\begin{equation*}
    J(G(u_j),u_j) \rightarrow m :=\inf_{u\in {H}_0^1(\Omega)}J(G(u),u) \text{ as } j \to \infty.
\end{equation*}
Then we aim to find a control $\overline{u}$ (with $\overline{y}=G(\overline{u})$) to which the minimizing sequence converges  such that $J(\overline{y},\overline{u}) = m$.

Since $\{u_j\}_{j\in \mathbb{N}}$ is a bounded sequence in $U\subset {H}_0^1(\Omega)$, there exists a weakly convergent subsequence $\{u_j\}_{j\in\mathcal{I}_1}$ with $\mathcal{I}_1\subset\mathbb{N} $ such that   
\begin{equation*}
    u_j \rightharpoonup \overline{u} \text{ in } {H}_0^1(\Omega) \text{ as }  j \to \infty\text{ in }\mathcal{I}_1 \text{ for some }\overline{u} \in {H}_0^1(\Omega).
\end{equation*} 
 As ${H}_0^1(\Omega)$ is compactly embedded into ${L}^2(\Omega)$, we apply the Rellich--Kondrachov theorem \cite[Chapter 5, Theorem 1]{evans} to get
\begin{equation*}
    u_j \rightarrow  \overline{u} \text{ in } {L}^2(\Omega) \text{ as }j \to \infty\text{ in }\mathcal{I}_1.
\end{equation*}
Since $\{u_j\}_{j\in\mathcal{I}_1} \subset U$ and  $\overline{u} \in {H}_0^1(\Omega)$, their ${H}_0^1(\Omega)$ norms are bounded by some constant $M$ such that we apply (\ref{equ:lip-g2}) to get
\begin{equation} \label{equ:tp1}
    \left\lVert \mathcal{G}(u_{j})-\mathcal{G}(\overline{u}) \right\rVert
    \leq 
    C(M)\left\lVert u_{j}-\overline{u} \right\rVert
    \to 0 , \text{ as } j \to \infty\text{ in }\mathcal{I}_1
\end{equation}
for some $C(M)>0$, which implies that $\{\mathcal{G}(u_{j})\}_{j\in\mathcal{I}_1}$ converges strongly in ${L}^2(\Omega)$.

By Theorem \ref{thm:wp-se}, $\|y_j\|_{{H}^2(\Omega)}\leq C(M)$ for $j\in \mathcal I_1$,  $\{y_{j}\}_{j\in\mathcal{I}_1}$ and thus $\{d(x,y_j)\}_{j\in\mathcal{I}_1}$ (based on \eqref{equ:bc-d}) are uniformly bounded in ${L}^\infty(\Omega)$ and hence in ${L}^2(\Omega)$. Therefore, there exists an $\mathcal{I}_2\subset\mathcal{I}_1 $ such that $d(x,y_{j}) \rightharpoonup \overline{d}$ in ${L}^2(\Omega)$ as $j\to\infty\text{ in }\mathcal{I}_2$. 
Furthermore,  (\ref{equ:tp1}) implies that $g(x,u_j) \to g(x,\overline{u})$ in ${L}^2(\Omega)$ as $j\to\infty\text{ in }\mathcal{I}_2$. Combine the above results to get $R_j :=g(x,u_j)-d(x,y_j) \rightharpoonup R:= g(x,\overline{u}) - \overline{d}$ in ${L}^2(\Omega)$ as $j\to\infty\text{ in }\mathcal{I}_2$, which, together with $y_j =\mathcal{A}^{-1} R_j$ and the fact that $\mathcal{A}^{-1} : {L}^2(\Omega) \to Y$ is a continuous linear operator \cite[Chapter 6, Theorem 1]{evans}, leads to
\begin{equation} \label{equ:y-wc}
        y_{j} \rightharpoonup \overline{y} \text{ in } Y \text{ as } j \to \infty\text{ in }\mathcal{I}_2\text{ for some } \bar y\in Y.
\end{equation}
By the compact embedding $Y \hookrightarrow \hookrightarrow  {L}^2(\Omega)$, we have the strong convergence in ${L}^2(\Omega)$
\begin{equation} \label{equ:y-sc}
        y_{j} \to  \overline{y} \text{ in } {L}^2(\Omega) \text{ as } j \to \infty\text{ in }\mathcal{I}_2.
\end{equation}

Based on (\ref{equ:lip-d}) in Assumption A, the Nemytskii operator $\mathcal{D}: y \to d(\cdot,y(\cdot))$ is locally Lipschitz continuous \cite[Lemma 4.11]{troltzschOptimalControlPartial2010}, that is, for any $E>0$ and $y_1,y_2$ satisfying $\left\lVert y_i \right\rVert_{{L}^\infty(\Omega)} \leq E$ ($i=1,2$), there exists a constant $C(E) > 0$ such that
\begin{equation*} 
    \left\lVert \mathcal{D}(y_1)-\mathcal{D}(y_2) \right\rVert 
    \leq 
    C(E) \left\lVert y_1-y_2 \right\rVert.
\end{equation*}
From (\ref{thm:wp-se}), the sequence $\{y_j\}_{j\in\mathcal{I}_2}$ and its limit $\overline{y}$ are uniformly bounded under the ${L}^\infty(\Omega)$ norm such that we invoke (\ref{equ:y-sc}) to get
\begin{equation} \label{equ:tp2}
   \left\lVert \mathcal{D}(y_j)-\mathcal{D}(\overline{y}) \right\rVert 
   \leq 
   C \left\lVert y_j-\overline{y} \right\rVert \to 0  \text{ as } j \to \infty\text{ in }\mathcal{I}_2.
\end{equation}

We then calculate the inner product of the state equation $\mathcal{A} y_j + d(x,y_j) = g(x,u_j)$ ($ j\in\mathcal{I}_2$) and $\phi\in Y$ to get
\begin{equation*}
    \int_{\Omega} a(x)\nabla y_j  \nabla\phi + c(x) y_j  \phi ~dx +  \int_{\Omega} d(x,y_j) \phi ~dx = \int_{\Omega} g(x,u_j) \phi ~dx.
\end{equation*}
Combine this with various convergence results, i.e. (\ref{equ:tp1})--(\ref{equ:y-wc}) and (\ref{equ:tp2}), to obtain 
\begin{equation*}
    \int_{\Omega} a(x)\nabla \overline{y} \, \nabla\phi + c(x) \overline{y}  \phi ~dx +  \int_{\Omega} d(x,\overline{y}) \phi ~dx = \int_{\Omega} g(x,\overline{u}) \phi ~dx.
\end{equation*}
Since $\phi$ is arbitrarily chosen in $Y$, we have $\overline{y}=G(\overline{u})$.

With the assistance of the above results, we turn to prove the main result. Let $J(y,u)=\Phi(y)+\Psi(u)$ where $\Phi(y) := \frac{1}{2}\left\lVert y-y_d \right\rVert^2 $ and $\Psi(u) := \frac{\lambda}{2} \left\lVert u \right\rVert_{{H}^1(\Omega)}^2$.

Since the ${H}^1_0$ norms of $\{u_j\}_{j\in\mathcal{I}_1} $ and  $\overline{u} $ are bounded by some constant $M$ (as mentioned above (\ref{equ:tp1})) and  $\{u_{j}\}_{j\in\mathcal{I}_1}$ converges strongly to $\bar u$ in ${L}^2(\Omega)$, 
 (\ref{equ:lip-G}) implies that $G(u_j) \to G(\overline{u}) = \overline{y}$ in ${L}^2(\Omega)$ as $j\to\infty$ in $\mathcal{I}_1$ such that
\begin{equation*}
    \lim_{j\to\infty} \Phi(G(u_j)) = \Phi(G(\overline{u})) = \Phi(\overline{y})\text{ as }j\to\infty\text{ in }\mathcal{I}_1.
\end{equation*}

Furthermore, as $\{u_{j}\}_{j\in\mathcal{I}_1}$ is a weakly convergent sequence to $ \overline{u}$  in ${H}_0^1(\Omega) \hookrightarrow {H}^1(\Omega)$, the convexity and weak lower semicontinuity of $\Psi$ imply 
\begin{equation*}
    \Psi(\overline{u}) \leq \liminf_{j \to \infty} \Psi(u_j).
\end{equation*}

Combine the above two estimates to get for $j\in\mathcal{I}_1$
\begin{equation*}
    \begin{aligned}
           m = \lim_{j\to\infty}J(G(u_j),u_j) =& \lim_{j\to\infty} \Phi(G(u_j)) +\liminf_{j \to \infty} \Psi(u_j) \ge  \Phi(\overline{y}) +  \Psi(\overline{u}) = J(\overline{y},\overline{u}),
    \end{aligned}
\end{equation*}
which, together with $m \leq J(\overline{y},\overline{u})$, implies the existence of an optimal control $\overline{u}$ with the associated state $\overline{y} = G(\overline{u})$.
\end{proof}


\section{Functional high-index saddle dynamics} \label{sec3}
Let $J(w)$ be a three times continuously Fr\'echet differentiable energy functional defined on $H_0^1(\Omega)$. For $w\in H_0^1(\Omega)$,  we write the
first, second, and third derivatives of $J(w)$ as $J'(w)$, $J''(w)$, and
$J'''(w)$, respectively. If we denote $\langle\cdot,\cdot\rangle$ as the duality pairing between
$H^{-1}(\Omega)$ and $H_0^1(\Omega)$, then we have $J'(w)[\phi]:=\langle J'(w),\phi\rangle$,
$J''(w)[\psi,\phi]:=\langle J''(w)\psi,\phi\rangle$ and $J'''(w)[\psi,\phi,\xi]:=\langle J'''(w)\psi\phi,\xi\rangle$ for $\phi,\psi,\xi\in H_0^1(\Omega)$. Then, by analogy with the traditional HiSD in Euclidean spaces \cite{yinHighIndexOptimizationBasedShrinking2019a}, we propose the following FHiSD that aims to locate an index-$k$ saddle point of $J$
\begin{equation}\label{equ:fhisd}
\left\{
\begin{aligned}
\dot{w} &= -J'(w) + 2\sum_{i=1}^{k} \langle J'(w), v_i \rangle v_i,\\
\dot{v}_i &= -J''(w)v_i + \langle J''(w)v_i, v_i \rangle v_i + 2\sum_{j=1}^{i-1} \langle J''(w)v_i, v_j \rangle v_j,\quad 1\le i\le k.
\end{aligned}
\right.
\end{equation}
Here $w$ is the primary variable and $\{v_i\}_{i=1}^k$ are directional variables used to identify the unstable subspace. Different from the traditional HiSD, where the gradient and Hessian of the energy function are explicitly given, the expressions of $J'(w),J''(w)v_i\in H^{-1}(\Omega)$ in (\ref{equ:fhisd}) may not be explicitly given, especially for optimal control problems as we will see later. In other words, the elements identifying these functionals based on the Riesz representation theorem are difficult to determine. Instead, the actions of the functionals $J'(w)$ and $J''(w)v_i$ on elements in $H^1_0(\Omega)$ can be evaluated (cf. Theorem \ref{thm:diff}), which again indicates the functional nature of the FHiSD. 

Throughout this section, we denote $z(x,t)$ as $z(t)$ for any function $z(x,t)$ and assume that $w(t),v_1(t),\ldots,v_k(t)\in H^1_0(\Omega)$ for each $t\geq 0$.
We first analyze the orthonormality preservation of the  FHiSD.

\begin{lemma}
If the initial values of (\ref{equ:fhisd}) satisfy $(v_i(0),v_j(0))=\delta_{ij}$ for $1\le i,j\le k$. Then $(v_i(t),v_j(t))=\delta_{ij}$ for $1\le i,j\le k$ and for $t\ge 0$.
\end{lemma}

\begin{proof}
We prove the result by induction on $i$. 
For $i=1$, we take the inner product of the equation for $v_1$ with $v_1$ to get
\begin{equation*}
\frac{d}{dt}\bigl((v_1,v_1)-1\bigr)
=2(\dot v_1,v_1)
=2\langle J''(w)v_1,v_1\rangle\bigl((v_1,v_1)-1\bigr).
\end{equation*}
Let $\rho(t):=(v_1,v_1)-1$ and $\mu(t):=2\langle J''(w)v_1,v_1\rangle$, then we have $\dot\rho(t)=\mu(t)\rho(t)$. Since $\rho(0)=0$, we obtain  $\rho(t)=0$ for all $t\ge 0$.

Assume that for some $2\leq m\leq k$, $(v_i(t),v_j(t))=\delta_{ij}$ for $1\le i,j\le m-1$ and all $t\ge 0$, and we intend to show $(v_m(t),v_j(t))=\delta_{mj}$ for $1\leq j\leq m$.
First, we prove $(v_m(t),v_j(t))=0$ for $1\le j\le m-1$ by induction on $j$. For $j=1$, taking inner products with $v_1$ and $v_m$ in the equations for $v_m$ and $v_1$, respectively, and using the symmetry of $J''(w)$, we obtain
\begin{equation*}
\frac{d}{dt}(v_m,v_1)
=(\dot v_m,v_1)+(\dot v_1,v_m)
=\bigl[\langle J''(w)v_m,v_m\rangle+\langle J''(w)v_1,v_1\rangle\bigr](v_m,v_1).
\end{equation*}
Since $(v_m(0),v_1(0))=0$, it follows that $(v_m(t),v_1(t))=0$ for all $t\ge 0$.

Next,  assume $(v_m(t),v_j(t))=0$ for $1\le j\le n-1$ for some $2\leq n\leq m-1$. Then
\begin{equation*}
\frac{d}{dt}(v_m,v_n)=(\dot v_m,v_n)+(\dot v_n,v_m)=\bigl[\langle J''(w)v_m,v_m\rangle+\langle J''(w)v_n,v_n\rangle\bigr](v_m,v_n).
\end{equation*}
Since $(v_m(0),v_n(0))=0$, we obtain $(v_m(t),v_n(t))=0$ for all $t\ge 0$. Thus, by inner induction, we have proved that $(v_m(t),v_j(t))=0$ for $1\le j\le m-1$.
The proof of $(v_m(t),v_m(t))=1$ can be performed in a similar manner as the proof of $(v_1(t),v_1(t))=1$. Thus we complete the proof.
\end{proof}

Now we prove the effectiveness of the FHiSD in locating saddle points by the following Theorem \ref{thm:ls-fhisd}, which extends the finite-dimensional counterpart for the traditional HiSD, cf. \cite[Theorem 1]{yinHighIndexOptimizationBasedShrinking2019a}. It is worth mentioning that in the finite-dimensional case, only the diagonal blocks of the Jacobian need to be analyzed, while for the infinite-dimensional case in Theorem \ref{thm:ls-fhisd}, the off-diagonal blocks need also to be bounded, which significantly complicates the analysis (cf. the detailed estimates in Theorem  \ref{thm:op}). Furthermore, instead of performing the matrix analysis in the  finite-dimensional case, a solution operator approach is adopted in the current  infinite-dimensional setting. For clarity, Theorem \ref{thm:ls-fhisd} is proved under general assumptions, which will be rigorously justified for the non-convex optimal control problems in Theorem  \ref{thm:op}.

The proof of Theorem \ref{thm:ls-fhisd} depends on the linearized system of (\ref{equ:fhisd}) at the stationary state $(w^*,v_1^*,\ldots,v_k^*)$, which takes the following form
\begin{equation}\label{lie}
\frac{d}{dt}
\begin{pmatrix}
\delta w\\
\delta v_1\\
\vdots\\
\delta v_k
\end{pmatrix}
=
\begin{pmatrix}
A      & 0      & \cdots & 0\\
B_1    & A_1    & \cdots & 0\\
\vdots & \vdots & \ddots & \vdots\\
B_k    & C_{k1} & \cdots & A_k
\end{pmatrix}
\begin{pmatrix}
\delta w\\
\delta v_1\\
\vdots\\
\delta v_k
\end{pmatrix}.
\end{equation}
Here, for $\phi\in H_0^1(\Omega)$ the block operators are given as below 
\begin{align}
A\phi &= -J''(w^*)\phi + 2\sum_{j=1}^{k}\langle J''(w^*)\phi,v_j^*\rangle v_j^*, \label{equ:ls-A-fhisd} \\
A_i\phi &= -J''\!(w^*)\phi + 2\langle J''\!(w^*)\phi,\!v_i^*\rangle v_i^* + \langle J''\!(w^*)v_i^*,\!v_i^*\rangle\phi + 2\sum_{j=1}^{i-1}\langle J''\!(w^*)\phi,\!v_j^*\rangle v_j^*, \label{equ:ls-Ai-fhisd} \\
B_i\phi &= -J'''(w^*)[\phi,v_i^*,\cdot] + J'''(w^*)[\phi,v_i^*,v_i^*]v_i^* + 2\sum_{j=1}^{i-1} J'''(w^*)[\phi,v_i^*,v_j^*]v_j^*, \label{equ:ls-Bi-fhisd}\\
C_{ij}\phi &= 2\langle J''(w^*)v_i^*,\phi\rangle v_j^* + 2\langle J''(w^*)v_i^*,v_j^*\rangle \phi,\quad 1\le j<i\le k.\label{equ:ls-Cij-fhisd}
\end{align}

\begin{theorem}\label{thm:ls-fhisd}
Let $w^*\in H_0^1(\Omega)$ and $\{v_i^*\}_{i=1}^k\subset H_0^1(\Omega)$ satisfy $(v_i^*,v_j^*)=\delta_{ij}$ for $1\leq i,j\leq k$. Suppose $(\gamma_1):$  there exist $\{\mu_m^*\}_{m\ge1}\subset\mathbb R\setminus\{0\}$ and an $L^2(\Omega)$-orthonormal basis $\{e_m^*\}_{m\ge1}\subset H_0^1(\Omega)$ of $L^2(\Omega)$ such that $\langle J''(w^*)e_m^*,\phi\rangle = \mu_m^*(e_m^*,\phi)$ for any $\phi\in H_0^1(\Omega)$ and $m\ge1$, and $(\gamma_2):$ if $\{v_i^*\}_{i=1}^k$ are orthonormal eigenfunctions of $J''(w^*)$, then the operators $B_i$ and $C_{ij}$ defined in \eqref{equ:ls-Bi-fhisd}--\eqref{equ:ls-Cij-fhisd} are bounded on $L^2(\Omega)$.
Then the following statements are equivalent:
\begin{itemize}
    \item[$(a)$] $(w^*,v_1^*,\ldots,v_k^*)$ is a stationary state of \eqref{equ:fhisd}, and the linearized system \eqref{lie} is exponentially stable in $L^2(\Omega)$;
    \item[$(b)$] $w^*$ is an index-$k$ saddle point of $J$, and $\{v_i^*\}_{i=1}^k$ are eigenfunctions associated with the negative eigenvalues of $J''(w^*)$.
\end{itemize}
\end{theorem}
\begin{proof}
We prove this theorem by two parts.
\vspace{4pt}

\noindent\textbf{Part I: Prove $\bm{(a)}$ from $\bm{(b)}$.}
Since $w^*$ is an index-$k$ saddle point, we have $J'(w^*)=0$ such that the right-hand side of the dynamics of $w$ in \eqref{equ:fhisd} vanishes at $w^*$. Moreover, there exist $\mu_1^*,\dots,\mu_k^*<0$ such that $\langle J''(w^*)v_i^*,\phi\rangle  = \mu_i^*(v_i^*,\phi)$ for $\forall \phi \in H_0^1(\Omega)$ and $1\le i\le k$. Using the orthonormality of $\{v_i^*\}_{i=1}^{k}$, we have $\langle J''(w^*)v_i^*,v_i^*\rangle=\mu_i^*$ and $\langle J''(w^*)v_i^*,v_j^*\rangle=0$ for $j<i$, which implies
\begin{equation*}
    -J''(w^*)v_i^*
    + \langle J''(w^*)v_i^*,v_i^*\rangle v_i^*
    + 2\sum_{j=1}^{i-1}\langle J''(w^*)v_i^*,v_j^*\rangle v_j^*
    = -J''(w^*)v_i^*+\mu_i^*v_i^*=0.  
\end{equation*}
Thus $(w^*,v_1^*,\ldots,v_k^*)$ is a stationary state of system \eqref{equ:fhisd}.

Next, we intend to prove the exponential stability of $(w^*,v_1^*,\ldots,v_k^*)$. Denote the other eigenvalues of $J''(w^*)$ by $\mu_{k+1}^*,\mu_{k+2}^*,\ldots$, and the corresponding eigenfunctions by $v_{k+1}^*,v_{k+2}^*,\ldots$. Then by assumption $(\gamma_1)$, $\{v_m^*\}_{m\ge1}$ form an orthonormal basis of $L^2(\Omega)$. 
 From \eqref{equ:ls-A-fhisd}, we immediately obtain
\begin{equation*}
    A v_m^*=\hat\mu_m^*v_m^*,\qquad \hat\mu_m^* =
\begin{cases}
\mu_m^*, & 1\le m\le k,\\
-\mu_m^*, & m\ge k+1.
\end{cases}
\end{equation*}
Since $w^*$ is an index-$k$ saddle point, we have $\mu_1^*<\cdots<\mu_k^*<0<\mu_{k+1}^*\le\cdots$. Hence $\hat\mu_m^*<0$ for all $m\ge1$. Let $\alpha_0:=-\sup_{m\ge1}\hat\mu_m^*>0$, and define the solution operator $S(t)$ associated with $A$ by $S(t)z:= \sum_{m=1}^{\infty}e^{\hat\mu_m^*t}(z,v_m^*)v_m^*$ for $z\in L^2(\Omega)$. 
By Parseval's identity,
\begin{equation*}
    \|S(t)z \|^2 = \sum_{m=1}^{\infty}e^{2\hat\mu_m^*t}|(z,v_m^*)|^2
    \le e^{-2\alpha_0 t}\sum_{m=1}^{\infty}|(z,v_m^*)|^2 = e^{-2\alpha_0 t}\|z\|^2.
\end{equation*}
Thus $S(t)$ is bounded on $L^2(\Omega)$ and satisfies $ \|S(t)z\|\le e^{-\alpha_0t}\|z\|$ for $
z\in L^2(\Omega)$ and $t\ge0$, which implies $\|\delta w\|\leq e^{-\alpha_0 t}\|\delta w(0)\|$.

To bound $\delta v_i$, we analyze the diagonal block $A_i$. 
Using \eqref{equ:ls-Ai-fhisd} and the orthonormality of $\{v_m^*\}_{m\ge1}$, we obtain
\begin{equation*}
   A_i v_m^*=\tilde\mu_{i,m}^*v_m^*,\qquad
   \tilde\mu_{i,m}^* = 
\begin{cases}
\mu_i^*+\mu_m^*, & 1\le m\le i,\\
\mu_i^*-\mu_m^*, & m\ge i+1.
\end{cases}
\end{equation*}
Since $w^*$ is an index-$k$ saddle point, we have $\tilde\mu_{i,m}^*<0$ for all $m\ge1$. Let $\alpha_i:=-\sup_{m\ge1}\tilde\mu_{i,m}^*>0$, and define the solution operators $S_i(t)$ associated with $A_i$ by $S_i(t)z := \sum_{m=1}^{\infty} e^{\tilde\mu_{i,m}^*t}(z,v_m^*)v_m^*$ for $z\in L^2(\Omega) $. Again by Parseval's identity, we have $ \|S_i(t)z\|\le e^{-\alpha_it}\|z\|$ for $ z\in L^2(\Omega)$ and $t\ge0$.
With $S_i(t)$ in hand, the equation of $\delta v_i$
\begin{equation*}
    \dot{\delta v_i}-A_i\delta v_i=B_i\delta w+\sum_{j=1}^{i-1}C_{ij}\delta v_j
\end{equation*}
can be reformulated as
\begin{equation}\label{jy1}
    \delta v_i(t)=S_i(t)\delta v_i(0)+\int_0^t S_i(t-s)\Big[B_i\delta w(s)+\sum_{j=1}^{i-1}C_{ij}\delta v_j(s)\Big]\,ds.
\end{equation}
Since $B_i$ and $C_{ij}$ are bounded on $L^2(\Omega)$ (cf. the assumption $(\gamma_2)$), there exist constants $M_i,M_{ij}>0$ depending only on $(w^*,v_1^*,\ldots,v_k^*)$ such that $\|B_i\xi\|\le M_i\|\xi\|$ and $\|C_{ij}\eta\|\le M_{ij}\|\eta\|$ for all $\xi,\eta\in H^1_0(\Omega)$. Let $\alpha_*:=\min\{\alpha_0,\alpha_1,\ldots,\alpha_k\}>0$, and we pick a $\beta$ satisfying $0<\beta<\alpha_*$. 

Now we bound $\delta v_i$ by induction on $i$. For $i=1$, we apply the boundedness of  $S_1(t)$ and $\delta w$ to get
\begin{equation*}
            \|\delta v_1(t)\| \le e^{-\alpha_1t}\|\delta v_1(0)\| +  M_1\int_0^t e^{-\alpha_1(t-s)} e^{-\alpha_0s}\|\delta w(0)\|\,ds.
\end{equation*}
Since $\beta<\min\{\alpha_0,\alpha_1\}$, there exists a constant $C>0$ such that $\int_0^t e^{-\alpha_1(t-s)} e^{-\alpha_0s}\,ds \le C e^{-\beta t}$. Thus $\|\delta v_1(t)\| \le K_1 e^{-\beta t} \bigl(\|\delta v_1(0)\|+\|\delta w(0)\|\bigr)$.

Assume that we have proved $\|\delta v_j(t)\| \le K_j e^{-\beta t} \Big( \|\delta w(0)\|+\sum_{\ell=1}^{j}\|\delta v_\ell(0)\| \Big)$ for $1\leq j\leq i-1$ for some positive constants $\{K_j\}_{j=1}^{i-1}$. Then we invoke these estimates in (\ref{jy1}) to get
\begin{equation*}
    \begin{aligned}
\|\delta v_i(t)\| &\le e^{-\alpha_i t}\|\delta v_i(0)\| + \int_0^t e^{-\alpha_i(t-s)} \Big( M_i\|\delta w(s)\| + \sum_{j=1}^{i-1}M_{ij}\|\delta v_j(s)\| \Big)ds  \\
&\le e^{-\alpha_i t}\|\delta v_i(0)\|
+ C\int_0^t e^{-\alpha_i(t-s)}e^{-\beta s}
\Big( \|\delta w(0)\|+\sum_{\ell=1}^{i-1}\|\delta v_\ell(0)\| \Big)ds.
\end{aligned}
\end{equation*}
Since $\beta<\alpha_i$, there exists a constant $C_i>0$ such that $\int_0^t e^{-\alpha_i(t-s)}e^{-\beta s}\,ds \le C_i e^{-\beta t}$. Therefore, $\|\delta v_i(t)\| \le K_i e^{-\beta t} \Big( \|\delta w(0)\|+\sum_{\ell=1}^{i}\|\delta v_\ell(0)\| \Big)$, which completes the induction. Consequently,
we have proved the exponential stability of the  linearized system (\ref{lie}) and thus the statement $(a)$.
\vspace{4pt}

\noindent\textbf{Part II: Prove $\bm{(b)}$ from $\bm{(a)}$.}
Define the projection $P^*z:=\sum_{i=1}^{k}\langle z,v_i^*\rangle v_i^*$ for $\forall z \in H^{-1}(\Omega)$. Since $\{v_i^*\}_{i=1}^{k}$ is orthonormal, $P^*$ is an orthogonal projection and $(I-2P^*)^2=I$. As $(w^*,v_1^*,\ldots,v_k^*)$ is a stationary state
of \eqref{equ:fhisd}, the equation of $w$ gives $ (I-2P^*) J'(w^*)=0$.
Applying $I-2P^*$ on both sides gives $J'(w^*)=0$. Thus $w^*$ is a stationary point of $J$.
Define $\mu_i^*:=\langle J''(w^*)v_i^*,v_i^*\rangle$ for $1\le i\le k$. Since $(w^*,v_1^*,\ldots,v_k^*)$ is a stationary state of \eqref{equ:fhisd}, we have $J''(w^*)v_i^* =\mu_i^*v_i^*+2\sum_{j=1}^{i-1}\langle J''(w^*)v_i^*,v_j^*\rangle v_j^*$. 
Taking the inner product with $v_\ell^*$, $1\le \ell<i$, and using the orthonormality of $\{v_i^*\}_{i=1}^k$, we obtain $\langle J''(w^*)v_i^*,v_\ell^*\rangle = 2\langle J''(w^*)v_i^*,v_\ell^*\rangle$, and hence $\langle J''(w^*)v_i^*,v_\ell^*\rangle=0$ for all $1\le \ell<i$. Therefore, $J''(w^*)v_i^*=\mu_i^*v_i^*$, that is,  $\langle J''(w^*)v_i^*,\phi\rangle  = \mu_i^*(v_i^*,\phi)$ for any $ \phi \in H_0^1(\Omega)$ and $1\le i\le k$.
The signs of the eigenvalues of $\{v_i^*\}_{i=1}^k$ can be determined from the exponential stability of the linearized system following the reverse procedure in Part I. Thus we omit the details of the derivation due to similarity.
\end{proof}

\section{Construction of control landscape}
\label{sec:hisd}
In this section, we justify the applicability of the FHiSD for non-convex optimal control, and then describe the numerical implementation of constructing the control landscape.
\subsection{FHiSD for optimal control}
By the control-to-state operator $G$, the cost functional $J$ could be reduced to $\hat{J}(u):=J(G(u),u)$. 
We first prove its differentiability and derive its derivatives, which will be directly used in the FHiSD.

\begin{theorem} \label{thm:diff}
Under the Assumptions A and B in the one-dimensional case (or Assumptions A and C in the two- and three-dimensional cases), the reduced cost functional $\hat{J}(u)$ is three times Fr\'echet differentiable on $H^1_0(\Omega)$, with derivatives defined as $\hat{J}'(u), \hat{J}''(u)$ and $\hat{J}'''(u)$. To characterize them, let $p \in Y$ denote the adjoint state uniquely determined by the adjoint equation associated with the state $y = G(u)$,
\begin{equation}\label{equ:ae-nc-ocp}
\mathcal A p+d_y(x,y)p=y-y_d
\text{ in }\Omega;~
p=0
\text{ on }\Gamma.
\end{equation}
Furthermore, let $z_i\in H_0^1(\Omega)$ for $1\leq i\leq 3$ and denote the first and second state derivatives as $\delta y_{z_i}:=G'(u)[z_i]$ $(1\leq i\leq 3)$ and $\delta y_{z_i z_j}:=G''(u)[z_i,z_j]$ $ (1\le i,j\le 3)$, which are determined by
\begin{align}\label{equ:linearized-state}
&\mathcal A \delta y_{z_i}+d_y(x,y) \delta y_{z_i}=g_u(x,u)z_i
\text{ in }\Omega;~
\delta y_{z_i} = 0
\text{ on }\Gamma,
\\
\label{equ:second-linearized-state}
&\mathcal A\delta y_{z_i z_j}+d_y(x,y)\delta y_{z_i z_j}= g_{uu}(x,u) z_i z_j-d_{yy}(x,y)\delta y_{z_i}\delta y_{z_j}
\text{ in }\Omega;~
\delta y_{z_i z_j} = 0
\text{ on }\Gamma.
\end{align}
Then, the derivatives $\hat{J}'(u), \hat{J}''(u)$ and $\hat{J}'''(u)$ are characterized by
\begin{align}\label{equ:df}
    &\langle \hat J'(u),z_1\rangle
    =
    \bigl(\lambda u+p\,g_u(x,u),z_1\bigr)
    +
    \lambda(\nabla u,\nabla z_1),\\
&    \langle \hat J''(u)z_1,z_2\rangle
    =
    \bigl(\lambda z_1+p\,g_{uu}(x,u) z_1, z_2\bigr)\nonumber\\
    &\hspace{0.9in}+
    \lambda(\nabla z_1,\nabla z_2)
    +
    \bigl((1-p\,d_{yy}(x,y))\delta y_{z_1}, \delta y_{z_2}\bigr),\label{equ:ddf}\\
&    \langle \hat J'''(u)z_1 z_2,z_3 \rangle 
    =
    \bigl((1-p\,d_{yy}(x,y))\delta y_{z_1 z_2}, \delta y_{z_3}\bigr) + \bigl((1-p\,d_{yy}(x,y))\delta y_{z_1 z_3},\delta y_{z_2}\bigr)\nonumber\\
    &\quad\quad\quad\quad\quad\quad\quad\quad+ 
    \bigl((1-p\,d_{yy}(x,y))\delta y_{z_2 z_3}, \delta y_{z_1}\bigr) + \bigl(p\,g_{uuu}(x,u)z_1 z_2 z_3,1\bigr)\nonumber\\
    &\quad\quad\quad\quad\quad\quad\quad\quad-
    \bigl(p\,d_{yyy}(x,y)\delta y_{z_1}\delta y_{z_2}\delta y_{z_3},1\bigr).\label{equ:dddf}
\end{align}

\end{theorem}

\begin{proof}
Fix $u_0 \in H_0^1(\Omega)$ and let $y_0:=G(u_0)\in Y:={H}^2(\Omega) \cap {H}^1_0(\Omega)$ be the associated state. Define $\mathcal S(y,u):=\mathcal Ay+d(x,y)-g(x,u)$ such that $\mathcal S(y_0,u_0)=0$. By the assumptions on $d$ and $g$, we have $\mathcal S \in C^3\bigl(Y \times H_0^1(\Omega),L^2(\Omega)\bigr)$. 
Then we consider $\mathcal S_y(y_0,u_0)\eta=\mathcal A\eta+d_y(x,y_0)\eta$ with $\eta\in Y$. 
Since $y_0\in Y\hookrightarrow L^\infty(\Omega)$, Assumption A implies $d_y(x,y_0)\in L^\infty(\Omega)$ and $d_y(x,y_0)\ge0$. Thus the bilinear form associated with $\mathcal A +d_y(x,y_0)$ is continuous and coercive on $H_0^1(\Omega)$. By the Lax--Milgram theorem \cite[Section 6.2, Theorem 1]{evans}, for every $f\in L^2(\Omega)$ there exists a unique weak solution $\eta\in H_0^1(\Omega)$ to $\mathcal A\eta+d_y(x,y_0)\eta=f$ with $\eta|_\Gamma=0$.
Furthermore, by elliptic regularity \cite[Section 6.3, Theorem 4]{evans}, we have $\eta \in Y$, and thus $\mathcal S_y(y_0,u_0):Y\to L^2(\Omega)$ is an isomorphism.
Combining this and $\mathcal S\in C^3\bigl(Y\times H_0^1(\Omega),L^2(\Omega)\bigr)$, the implicit function theorem \cite[Appendix C, Theorem 9]{evans} implies that there exists a neighborhood $U_0\subset H_0^1(\Omega)$ of $u_0$ and a unique mapping $\widetilde G\in C^3(U_0,Y)$ such that $\mathcal S(\widetilde G(u),u)=0$ for any $u\in U_0$.
By the uniqueness of the state equation, $\widetilde G$ coincides with the original control-to-state mapping $G$ on $U_0$. Hence $G\in C^3(U_0,Y)$, which, together with the chain rule, leads to $\hat J\in C^3(U_0)$.
Since $u_0$ is arbitrary, $\hat J$ is three times Fr\'echet differentiable on $H_0^1(\Omega)$.
The detailed derivation of \eqref{equ:df}--\eqref{equ:dddf} is given in Appendix~\ref{app:derivative-proof}.
\end{proof}

Based on Theorem \ref{thm:diff}, we have $\hat J'(u)\in H^{-1}(\Omega)$ and $\hat J''(u):H_0^1(\Omega)\to H^{-1}(\Omega)$ such that the FHiSD for the optimal control problem $\hat J(u)$ reads:
\begin{equation}\label{equ:hisd-ocp}
\begin{aligned}
    &\text{FHiSD for}\\
    &\text{~optimal}\\
    &\text{\,\,\,control}
\end{aligned}
\left\{
\begin{aligned}
    &\text{FHiSD}
    \left\{
    \begin{aligned}
        \!\dot u &\!=\! -\hat J'(u) +  2\sum_{i=1}^{k}\langle \hat J'(u),v_i\rangle v_i, \\
        \!\dot v_i &\!=\! -\hat J''\!(u)v_i \!+\! \langle \hat J''\!(u)v_i,\!v_i\rangle v_i \!+ \!2\sum_{j=1}^{i-1}\langle \hat J''\!(u)v_i,\!v_j\rangle v_j,~~~ \!1\le i\le k.
        \end{aligned}
        \right.\\
    &\begin{aligned}
    &\text{Functional}\\
    &\text{evaluation}
     \end{aligned}
    \left\{
        \begin{aligned}
&\langle \hat{J}'(u),v_i\rangle
    = \bigl(\lambda u + p\, g_u(x,u), v_i\bigr) + \lambda (\nabla u, \nabla v_i),\\
&   \begin{aligned}
    \langle \hat{J}''(u)v_i,v_j\rangle
    &= \bigl(\lambda v_i + p\, g_{uu}(x,u)v_i,v_j\bigr)\\
    &\quad \!+ \lambda (\nabla v_i,\!\nabla v_j)  +  \bigl( (1 \!-\! p\, d_{yy}(x,y)) \delta y_{v_i}, \delta y_{v_j} \bigr).
    \end{aligned}
    \end{aligned}
    \right.\\
    &\begin{aligned}
    &\text{Constraint}\\
    &\text{~~~PDEs}
     \end{aligned}
    \left\{
        \begin{aligned}
&\mathcal{A} y + d(x,y) = g(x,u) \text{ in } \Omega,\text{ with } y = 0 \text{ on } \Gamma.\\
&\mathcal{A} p + d_y(x,y)p = y - y_d \text{ in } \Omega, \text{ with } p = 0 \text{ on } \Gamma.\\
&\mathcal{A} \delta y_{v_i} + d_y(x,y)\delta y_{v_i} = g_u(x,u)v_i \text{ in } \Omega,\text{ with } \delta y_{v_i} \!=\! 0 \text{ on } \Gamma,\\
&1\le i\le k.
        \end{aligned}
    \right.
\end{aligned}
\right.
\end{equation}

\subsection{Effectiveness of FHiSD for optimal control}

To validate the effectiveness of (\ref{equ:hisd-ocp}) in locating stationary points of the non-convex optimal control problem (\ref{equ:nc-ocp})--(\ref{equ:se-nc-ocp}), we need to justify the assumptions $(\gamma_1)$--$(\gamma_2)$ in Theorem~\ref{thm:ls-fhisd} for the reduced functional $\hat J$. 
We first collect several auxiliary estimates that will be used to verify the assumptions of Theorem~\ref{thm:ls-fhisd}.

\begin{lemma}\label{lm:ae}
Let $u^*\in H_0^1(\Omega)$ and $y^*=G(u^*)$. Let $p^*$ be the solution of the adjoint equation \eqref{equ:ae-nc-ocp}, and set $\mathcal L:=\mathcal A+d_y(x,y^*)$. Suppose the Assumptions A and B hold in the one-dimensional case or Assumptions A and C hold in the two- and three-dimensional cases. Then the following assertions hold:

\begin{itemize}
\item[$(a)$] The coefficients $d_y(x,y^*)$, $d_{yy}(x,y^*)$, $d_{yyy}(x,y^*)$, $g_u(x,u^*)$, $g_{uu}(x,u^*)$ and $g_{uuu}(x,u^*)$ are bounded in $L^\infty(\Omega)$.

\item[$(b)$] The $p^*\in Y\hookrightarrow L^\infty(\Omega)$ and $\|p^*\|_{H^2(\Omega)}\le C\|y^*-y_d\|$. 

\item[$(c)$] For $\eta\in L^2(\Omega)$, the linearized equation 
\begin{equation*}
        \mathcal L\delta y_\eta=g_u(x,u^*)\eta \text{ in }\Omega;~
        \delta y_\eta=0 \text{ on }\Gamma 
\end{equation*} 
admits a unique solution $\delta y_\eta\in Y$ satisfying $\|\delta y_\eta\|_{H^2(\Omega)}\le C\|\eta\|$.

\item[$(d)$] For  $v\in H_0^1(\Omega)\cap L^\infty(\Omega)$ and $\eta\in L^2(\Omega)$, the second-order linearized equation 
\begin{equation*}
    \mathcal L\delta y_{\eta v} = g_{uu}(x,u^*)\eta v - d_{yy}(x,y^*)\delta y_\eta\delta y_v \text{ in }\Omega;~
    \delta y_{\eta v}=0 \text{ on }\Gamma 
\end{equation*}
admits a unique solution $\delta y_{\eta v}\in Y$ satisfying   $\|\delta y_{\eta v}\|_{H^2(\Omega)} \le C\|\eta\|$, where $C>0$ is independent of $\eta$.
\end{itemize}
\end{lemma}

\begin{proof}
We first prove $(a)$. We combine $y^*\in Y \hookrightarrow L^\infty(\Omega)$ for $n\le3$ and the Assumption A to find that the derivatives of $d$ in $(a)$ are bounded.
The boundedness of the derivatives of $g$ in $(a)$ in the one-dimensional case follows from the Assumption B and $u^*\in H_0^1(\Omega)\hookrightarrow L^\infty(\Omega)$.
In the two- and three-dimensional cases, Assumption C with $j=0$ implies $|g_s(x,s)|\le L_0$ for any $ s\in\mathbb R$ and a.e. $x\in\Omega$.
In particular, we have $|g_u(x,u^*(x))|\le L_0~\text{for a.e. }x\in\Omega$, that is, $g_u(x,u^*)\in L^\infty(\Omega)$. The boundedness of $g_{uu}(x,u^*)$ and $g_{uuu}(x,u^*)$ can be proved similarly.

By $(a)$ and Assumption A, $d_y(x,y^*)$ is bounded and nonnegative. Hence the operator $\mathcal L$ is coercive such that for any $f\in L^2(\Omega)$ the equation $\mathcal L z=f$ with $ z|_{\Gamma}=0$ admits a unique solution $z\in Y$ satisfying $\|z\|_{H^2(\Omega)}\le C\|f\|$.
Applying this estimate, we can prove $(b)$--$(d)$. Specifically, the conclusion $(b)$ is a direct consequence of this estimate. To prove $(c)$, we note that the boundedness of $g_u(x,u^*)$ implies $\|g_u(x,u^*)\eta\|\le C\|\eta\|$ for $\eta\in L^2(\Omega)$, leading to $\|\delta y_\eta\|_{H^2(\Omega)} \le C\|g_u(x,u^*)\eta\| \le C\|\eta\|$.
This proves $(c)$.
To prove $(d)$, let $v\in H_0^1(\Omega)\cap L^\infty(\Omega)$ and $\eta\in L^2(\Omega)$. By $(c)$, $\delta y_\eta\in Y$ and $\|\delta y_\eta\|_{H^2(\Omega)}\le C\|\eta\|$. Also, since $v\in L^2(\Omega)$, $\delta y_v\in Y\hookrightarrow L^\infty(\Omega)$.
We combine these with the boundedness of $g_{uu}(x,u^*)$ and $d_{yy}(x,y^*)$ to get
\begin{equation*}
    \|g_{uu}(x,u^*)\eta v
    -d_{yy}(x,y^*)\delta y_\eta\delta y_v\|
    \le
    C\|v\|_{L^\infty(\Omega)}\|\eta\|
    +
    C\|\delta y_v\|_{L^\infty(\Omega)}\|\delta y_\eta\| 
    \le C\|\eta\|,
\end{equation*}
which implies $\|\delta y_{\eta v}\|_{H^2(\Omega)} \le C\|\eta\|$.
This proves $(d)$.
\end{proof}

\begin{theorem}[Justification of assumptions $(\gamma_1)$--$(\gamma_2)$ in Theorem~\ref{thm:ls-fhisd} for $\hat J$]
\label{thm:op}
   Let $u^*\in H_0^1(\Omega)$ and $\{v_i^*\}_{i=1}^k\subset H_0^1(\Omega)$ satisfy $(v_i^*,v_j^*)=\delta_{ij}$ for $1\leq i,j\leq k$.
Suppose the Assumptions A and B hold in the one-dimensional case or Assumptions A and C hold in the two- and three-dimensional cases. 
Then if $\hat J''(u^*)$ is non-degenerate (i.e. $\hat J''(u^*)$ does not have zero eigenvalues), the assumption $(\gamma_1)$ of Theorem~\ref{thm:ls-fhisd} is valid for $J=\hat J$.
If in addition $\{v_i^*\}_{i=1}^k\subset H_0^1(\Omega)\cap L^\infty(\Omega)$, then the assumption $(\gamma_2)$ is valid for $J=\hat J$.
\end{theorem}

\begin{proof}
For $u^*\in H^1_0(\Omega)$, let $y^*=G(u^*)$ and $p^*$ be the corresponding state variable and adjoint state, respectively. Then we apply Theorem \ref{thm:wp-se} and Lemma~\ref{lm:ae}$(b)$ to obtain $y^*,p^*\in Y \hookrightarrow L^\infty(\Omega)$.
Let $q^*:=1-p^*d_{yy}(x,y^*)$, which satisfies $q^*\in L^\infty(\Omega)$ due to $d_{yy}(x,y^*) \in L^\infty(\Omega)$.
Then we justify the assumptions $(\gamma_1)$--$(\gamma_2)$ in Theorem~\ref{thm:ls-fhisd} in the following two parts.
\vspace{4pt}

\noindent\textbf{Part I: Validation of assumption $(\gamma_1)$.}
Define the symmetric bilinear form $a_*:H_0^1(\Omega)\times H_0^1(\Omega)\to \mathbb R$ by $a_*(v,\phi):=\langle \hat J''(u^*)v,\phi\rangle$. For a sufficiently large constant $\rho>0$, we introduce the shifted form $a_\rho(v,\phi) :=a_*(v,\phi)+\rho(v,\phi)$. By Theorem~\ref{thm:diff}, for any $v,\phi\in H_0^1(\Omega)$, the shifted form admits the representation
\begin{equation*}
        a_\rho(v,\phi) =  \lambda\big[(\nabla v,\nabla\phi) + (v,\phi)\big] + (p^*g_{uu}(x,u^*)v,\phi) + \bigl(q^*\delta y_v,\delta y_\phi\bigr) + \rho(v,\phi).
\end{equation*}
The first right-hand side term $\lambda\big[(\nabla v,\nabla\phi) +(v,\phi)\big]$ is symmetric, continuous, and coercive on $H_0^1(\Omega)$. 
Moreover, by Lemma~\ref{lm:ae}$(a)$, we have $g_{uu}(x,u^*)\in L^\infty(\Omega)$; thus the second right-hand side term is symmetric, and there exists $C_1>0$ such that $|(p^*g_{uu}(x,u^*)v,\phi)| \le C_1\|v\|\|\phi\|$.
For the third right-hand side term, since $v,\phi \in H_0^1(\Omega) \! \subset \! L^2(\Omega)$, Lemma~\ref{lm:ae}$(c)$ gives $\|\delta y_v\|_{H^2(\Omega)}\le C\|v\|$, and $\|\delta y_\phi\|_{H^2(\Omega)}\le C\|\phi\|$.  
Hence, there exists $C_2>0$ such that 
\begin{equation*}
    \bigl|(q^*\delta y_v,\delta y_\phi)\bigr|
    \le \|q^*\|_{L^\infty(\Omega)} \|\delta y_v\| \|\delta y_\phi\|
    \le C_2\|v\| \|\phi\|.
\end{equation*}

With these estimates, let $C_0=C_1+C_2$, and choose $\rho>0$ sufficiently large such that
$\lambda-C_0+\rho>0$, and set $c:=\min\{\lambda,\lambda - C_0 + \rho\}>0$. Then we take $\phi=v\in H_0^1(\Omega)$ in $a_\rho(v,\phi)$ to get the coercivity
\begin{equation*}
        a_\rho(v,v) \ge \lambda \|\nabla v\|^2 + (\lambda - C_0 + \rho)\|v\|^2  \ge c\bigl(\|\nabla v\|^2+\|v\|^2\bigr) = c \|v\|^2_{H^1(\Omega)} .
\end{equation*}
The continuity can be proved similarly. Thus by the Lax--Milgram theorem, for each $f\in L^2(\Omega)$ there exists a unique $z\in H_0^1(\Omega)$ satisfying $a_\rho(z,\phi)=(f,\phi)$ for any $\phi\in H_0^1(\Omega)$. Now we define the solution operator $T:L^2(\Omega)\to H_0^1(\Omega)$ by $Tf=z$.
By the symmetry of $a_\rho$, the operator $T$ is self-adjoint. Since the embedding $H_0^1(\Omega)\hookrightarrow L^2(\Omega)$ is compact, $T$ is compact on $L^2(\Omega)$. Therefore, by the spectral theorem for compact self-adjoint operators \cite[Appendix D, Theorem 7]{evans}, there exist a countable $L^2(\Omega)$-orthonormal basis $\{e_m^*\}_{m\ge1}$ and positive eigenvalues $\{\nu_m\}_{m\ge1}$ such that $T e_m^* = \nu_m e_m^*$. By the definition of $T$ and $a_\rho(e_m^*,\phi) = a_*(e_m^*,\phi) + \rho(e_m^*,\phi)$, we obtain 
\begin{equation*}
    a_*(e_m^*,\phi) := \langle \hat J''(u^*)e_m^*,\phi\rangle = \mu_m^*(e_m^*,\phi),\quad \forall \phi\in H_0^1(\Omega),
\end{equation*}
where $\mu_m^*:=\nu_m^{-1}-\rho$, which, together with the non-degeneracy of $\hat J''(u^*)$, justifies the assumption $(\gamma_1)$ in Theorem~\ref{thm:ls-fhisd}.

\vspace{6pt}

\noindent\textbf{Part II: Validation of assumption $(\gamma_2)$.}
We first bound $C_{ij}$. Since $\langle \hat J''(u^*)v_i^*,\phi\rangle  = \hat\mu_i^*(v_i^*,\phi)$ for any $\phi \in H_0^1(\Omega)$ and for some $\hat\mu_i^*\in\mathbb R\setminus\{0\}$ as presumed in assumption $(\gamma_2)$, 
we have $|\langle \hat J''(u^*)v_i^*,\phi\rangle|\le |\hat\mu_i^*|\,\|v_i^*\|\,\|\phi\|$ for any $  \phi\in H_0^1(\Omega)$.
Using the orthonormality one has $C_{ij}\eta = 2\hat\mu_i^*(v_i^*,\eta)v_j^*$ for $j<i$ and for any $\eta\in H^1_0(\Omega)$, leading to
$\|C_{ij}\eta\| \le  2|\hat\mu_i^*|\,\|\eta\|$. Hence $C_{ij}$ is bounded on $L^2(\Omega)$.

We turn to bound $B_i$. By definition in \eqref{equ:ls-Bi-fhisd}, for any $\xi\in H_0^1(\Omega)$ and $\{v_i^*\}_{i=1}^k \subset H_0^1(\Omega)\cap L^\infty(\Omega)$, we have
\begin{equation*}
B_i\xi = -\hat J'''(u^*)[\xi,v_i^*,\cdot] +\hat J'''(u^*)[\xi,v_i^*,v_i^*]v_i^*  +
2\sum_{j=1}^{i-1}\hat J'''(u^*)[\xi,v_i^*,v_j^*]v_j^* .
\end{equation*}

Let $\zeta\in L^2(\Omega)$. Following \eqref{equ:dddf}, we have
\begin{equation*}
\begin{aligned}
\hat J'''(u^*)[\xi,v_i^*,\zeta]
&=
\bigl(q^*\delta y_{\xi v_i^*},\delta y_\zeta\bigr) 
+\bigl(q^*\delta y_{\zeta v_i^*},\delta y_\xi\bigr)
+\bigl(q^*\delta y_{\xi\zeta},\delta y_{v_i^*}\bigr) \\
&\quad+ \bigl(p^*g_{uuu}(x,u^*)\xi v_i^*\zeta,1\bigr)
-\bigl(p^*d_{yyy}(x,y^*)\delta y_\xi\delta y_{v_i^*}\delta y_\zeta,1\bigr)\\
&=: I_1 + I_2 + I_3 + I_4 - I_5.
\end{aligned}
\end{equation*}

For $I_1$ and $I_2$, since $\xi,\zeta\in L^2(\Omega)$ and $v_i^* \in H_0^1(\Omega)\cap L^\infty(\Omega)$, Lemma~\ref{lm:ae} $(c)$--$(d)$ imply $\|\delta y_\xi\|_{H^2(\Omega)}\le C\|\xi\|$ and $\|\delta y_{\xi v_i^*}\|_{H^2(\Omega)}\le C_i\|\xi\|$. 
Therefore, we have $|I_1| \le \|q^*\|_{L^\infty(\Omega)}\|\delta y_{\xi v_i^*}\|\|\delta y_\zeta\| \le C_i\|\xi\|\,\|\zeta\|$.
Similarly, $|I_2|\le C_i\|\xi\|\,\|\zeta\|$.
For $I_3$, since a direct estimate of $\delta y_{\xi\zeta}$ by Lemma~\ref{lm:ae} (d) would require stronger regularity of $\zeta$, we adopt an alternative approach. Let $r \in Y$ solve $\mathcal L r =q^*\delta y_{v_i^*}$ with $ r|_{\Gamma}=0$, and thus
    $I_3 = \bigl(\delta y_{\xi\zeta},\mathcal L r\bigr)
    = \bigl(\mathcal L\delta y_{\xi\zeta},r\bigr) 
    = \bigl(g_{uu}(x,u^*)r\xi,\zeta\bigr)
    - \bigl(d_{yy}(x,y^*)\delta y_\xi r,\delta y_\zeta\bigr).$
Hence, using the boundedness of $g_{uu}(x,u^*)$, $d_{yy}(x,y^*)$, and $r$, we get
    $|I_3|\!
    \le C\|\xi\|\,\|\zeta\|
    + C\|\delta y_\xi\|\,\|\delta y_\zeta\|
    \le C\|\xi\|\,\|\zeta\|$.
    For $I_4$ and $I_5$, since $g_{uuu}(x,u^*)$ and $d_{yyy}(x,y^*)$ are bounded, we have $|I_4| \le C\|\xi\|\,\|\zeta\|$ and $|I_5|\le C\|\delta y_\xi\|\,\|\delta y_\zeta\|\le C\|\xi\|\,\|\zeta\|$.
Combining the estimates for all five terms yields $|\hat J'''(u^*)[\xi,v_i^*,\zeta]| \le C\|\xi\|\,\|\zeta\|$ for any $\zeta\in L^2(\Omega)$, which implies $\|\hat J'''(u^*)[\xi,v_i^*,\cdot]\|  \le C\|\xi\|$.
We combine this with $\|v_i^*\|=1$ to obtain $
    \|B_i\xi\|\le M_i\|\xi\|$ for some $M_i>0$, which completes the proof.
\end{proof}


\subsection{Numerical implementation of control landscape}
We first discretize the FHiSD. Let $\mathcal T_h$ be a  partition of $\Omega$ with mesh size $h$, and define
\vspace{-2pt}
\begin{equation*}
V_h:=\bigl\{\phi_h\in C(\overline{\Omega})\cap H^1_0(\Omega) : \phi_h|_K\in P_1(K),\, \forall K\in\mathcal T_h\bigr\},
\end{equation*}
where $P_1(K)$ denotes the space of linear functions on $K$. 
We denote the discrete control, state, adjoint and directional variables at the $n$-th iteration by $u_h^{(n)}$, $y_h^{(n)}$, $p_h^{(n)}$, and $v_{h,i}^{(n)}$ ($1\le i\le k$), respectively.
Let $a(\cdot,\cdot)$ denote the bilinear form associated with the elliptic operator $\mathcal A$. 
Then the finite element scheme of the FHiSD for optimal control (\ref{equ:hisd-ocp}) at the $n$-th iteration reads:
\vspace{-2pt}
\begin{equation}\label{equ:dhisd-ocp}
\begin{aligned}
    &\text{Discrete}\\
    &\text{~FHiSD}\\
    &\text{~\,~~for}\\
    &\text{optimal}\\
    &\text{\,control}
\end{aligned}
\left\{
\begin{aligned}
    &\!\begin{aligned}
    &\text{FHiSD}\\
    &\text{scheme}
    \end{aligned}
    \!\left\{
        \begin{aligned}
&(u_h^{(n+1)}\!,\!\phi_h) \!= (u_h^{(n)},\phi_h)\\
&\quad \quad \quad \quad \quad \quad -\beta_n \Big\{\!\langle \hat J'(u_h^{(n)}),\!\phi_h\rangle \!-\! 2 \sum_{i=1}^{k} \langle\hat J'(u_h^{(n)}),\! {v}^{(n)}_{h,i}\rangle ({v}^{(n)}_{h,i},\!\phi_h)\!\Big\}\!,\\
&({\hat{v}}^{(n+1)}_{h,i},\!\psi_h) \!=  ({v}^{(n)}_{h,i},\psi_h)\\
&\quad \quad \quad \quad \quad \quad  - \gamma_n \langle \hat J''(u_h^{(n)})v_{h,i}^{(n)},v_{h,i}^{(n)}\rangle({v}^{(n)}_{h,i},\psi_h),~~~~~~~~1\le i\le k, \\
&\{{v}^{(n+1)}_{h,i}\} _{i=1}^k= \text{ Orth } \{{\hat{v}}^{(n+1)}_{h,i}\} _{i=1}^k.
\end{aligned}
\right.\\
&\!\begin{aligned}
    &\text{Functional}\\
    &\text{Evaluation}
     \end{aligned}
    \!\left\{
        \begin{aligned}
&\!\langle \hat J'\!(u_h^{(n)}),\phi_h\rangle\! =\! (p_h^{(n)}\!\!g_u(x,\!u_h^{(n)}),\phi_h \!) \!+\! \lambda( u_h^{(n)}\!,\!\phi_h \!)\!+\!\lambda(\nabla \!u_h^{(n)}\! ,\!\nabla\! \phi_h\!),\\
&\begin{aligned}
        \!\langle \hat J''\!(u_h^{(n)})v_{h,i}^{(n)},\psi_h\rangle & = (p_h^{(n)}\,g_{uu}(x,u_h^{(n)})v_{h,i}^{(n)}, \psi_h) \\
        &\quad +\lambda (v_{h,i}^{(n)},\psi_h) + \lambda (\nabla v_{h,i}^{(n)},\nabla \psi_h)  \\
        &\quad +((1-p_h^{(n)}\,d_{yy}(x,y_h^{(n)})) \delta y_{h,i}^{(n)},\delta y_{\psi_h}^{(n)}).
    \end{aligned}
    \end{aligned}
    \right.\\
&\!\begin{aligned}
    &\text{~Constraint}\\
    &\text{PDE solvers}
    \end{aligned}
    \!\left\{
        \begin{aligned}
&a(y_h^{(n)},\chi_h) + (d(y_h^{(n)}),\chi_h) = (g(u_h^{(n)}),\chi_h),\\
&a(p_h^{(n)},\chi_h) + (d_y(x,y_h^{(n)}) p_h^{(n)},\chi_h) = (y_h^{(n)} - y_d,\chi_h),\\
&a(\delta y_{h,i}^{(n)},\chi_h) + (d_y(x,y_h^{(n)}) \delta y_{h,i}^{(n)},\chi_h) = (g_u(x,u_h^{(n)})v_{h,i}^{(n)},\chi_h),\\
&a(\delta y_{\psi_h}^{(n)},\chi_h) + (d_y(x,y_h^{(n)}) \delta y_{\psi_h}^{(n)},\chi_h) = (g_u(x,u_h^{(n)})\psi_h,\chi_h).
        \end{aligned}
    \right.\\
\end{aligned}
\right.
\end{equation}
Here $\beta_n$ and $\gamma_n$ denote the time-step sizes, and the operator ``Orth'' refers to the standard Gram--Schmidt orthonormalization process. It is worth mentioning that although the discrete schemes for the directional variables in \eqref{equ:dhisd-ocp} drop nonlinear terms about $v_i$ and $v_j$ in their continuous formulations in \eqref{equ:fhisd}, these terms are indeed involved in the Gram--Schmidt process, as shown in \cite{Z3CMS}. In addition, we adopt the Barzilai--Borwein (BB) method for adaptive step-size selection \cite{barzilaiTwoPointStepSize1988,yinHighIndexOptimizationBasedShrinking2019a}, and apply Newton's method to solve the nonlinear state equation.


With the assistance of the FHiSD in locating  saddle points of the optimal control problem (\ref{equ:nc-ocp})--(\ref{equ:se-nc-ocp}), we explain how to construct its control landscape. Following \cite{yinSearchingSolutionLandscape2021b}, the downward and upward search algorithms will be adopted, where the former starts from a high-index saddle point to locate saddle points of lower indices and the latter ascends from a low-index point toward saddle points of higher indices. By combining these two search methods, we could systematically locate multiple saddle points of the optimal control problem to generate its control landscape. 

We present details for the downward search, which is frequently used in constructing the control landscape from a parent state, as we will do in numerical experiments (more details on downward and upward searches can be found in \cite{yinSearchingSolutionLandscape2021b}).   Suppose that $u_h^{\left\langle k\right\rangle}$ represents an index-$k$ saddle point with eigenvectors (stored column-wise in  $V_h^{\left\langle k\right\rangle}:=[v_{h,1},v_{h,2},...,v_{h,k}]$) constituting unstable subspaces.
To search for an index-$r$ saddle point ($r<k$) from $u_h^{\left\langle k\right\rangle}$,
we first select one direction $v_{h,i}$ with $1\leq i\leq r$ from $V_h^{\left\langle k\right\rangle}$ and perturb
$u_h^{\left\langle k\right\rangle}$ along $v_{h,i}$ to obtain the initial state of control variable $u_h^{\left\langle r\right\rangle}$, that is, $u_h^{\left\langle r\right\rangle}=u_h^{\left\langle k\right\rangle}\pm \sigma v_{h,i}$ for some $\sigma>0$.
Then we choose $r$ unstable directions  from the remaining set $V_h^{\left\langle k\right\rangle}\setminus \{v_{h,i}\}$ to get the initial directions $V_h^{\left\langle r\right\rangle }:=\{\tilde{v}_{h,1},\tilde{v}_{h,2},...,\tilde{v}_{h,r}\}$.
With these initial conditions, we can apply \eqref{equ:dhisd-ocp} to locate an index-$r$ saddle point.

\section{Numerical experiments}
\label{sec:experiments}
We present several numerical experiments and comparisons with classical methods to demonstrate the advantages of the control landscape method. When the index of the parent state is high and there are many intermediate saddle points, we omit most of them and focus on low-index ones in the control landscape. 

\subsection{Numerical comparison for one-dimensional problem}
Let $\Omega=[0,1]$, $a_{11}(x)=c(x)\equiv 1$, $ d(x,y) = y^3$, $g(x,u) = 0.001u^2 + \cos(2\pi u)$ and $y_d = -2\sin(\pi x)$. A uniform partition of $\Omega$ is adopted. Before performing the numerical comparison, we illustrate the mesh independence of the numerical results to determine a suitable mesh size $h$. Let  $\lambda=0.02$, and we present the stationary points in control landscapes of (\ref{equ:nc-ocp})--(\ref{equ:se-nc-ocp}) computed under different $h$ in Table \ref{tab:md}. We find that under successive mesh refinements, the Morse indices of the stationary points remain unchanged, the residuals measured by $\|\hat J'\| $ at these  points remain stable, and the functional values $\hat J$ at these points converge. These results suggest the mesh independence of the numerical results for small $h$, and we pick $h=2^{-8}$ for one-dimensional problems. 

\begin{table}[htbp]
\centering
\setlength{\tabcolsep}{5pt}
\setlength{\belowcaptionskip}{2pt}
{\footnotesize
\caption{Mesh independence test of the Morse indices of stationary points and the corresponding functional values and residuals.}
\label{tab:md}
\begin{tabular}{ccc|ccc|ccc|ccc}
\hline
\multicolumn{3}{c|}{\rule{0pt}{8pt} $h=2^{-5}$} & \multicolumn{3}{c|}{$h=2^{-6}$}  & \multicolumn{3}{c|}{$h=2^{-7}$} & \multicolumn{3}{c}{$h=2^{-8}$} \\
\hline
  \rule{0pt}{8pt}
  $k$ & $\hat J$ & $\|\hat J'\|$
& $k$ & $\hat J$ & $\|\hat J'\|$ 
& $k$ & $\hat J$ & $\|\hat J'\|$  
& $k$ & $\hat J$ & $\|\hat J'\|$ \\
\hline
4 & 1.1187 & 0        & 4 & 1.1200 & 0        & 4 & 1.1203 & 0        & 4 & 1.1204 & 0 \\
3 & 1.1046 & 9.257E-5 & 3 & 1.1051 & 9.859E-5 & 3 & 1.1053 & 9.137E-5 & 3 & 1.1053 & 9.992E-5 \\
2 & 1.0598 & 9.895E-5 & 2 & 1.0603 & 9.460E-5 & 2 & 1.0605 & 9.516E-5 & 2 & 1.0605 & 9.729E-5 \\
1 & 0.9955 & 9.645E-5 & 1 & 0.9963 & 9.885E-5 & 1 & 0.9965 & 9.641E-5 & 1 & 0.9966 & 9.798E-5 \\
0 & 0.9209 & 9.339E-5 & 0 & 0.9219 & 9.437E-5 & 0 & 0.9222 & 9.864E-5 & 0 & 0.9222 & 9.801E-5 \\
\hline
\end{tabular}}
\end{table}

\begin{figure}[h]
    \centering
    \includegraphics[width=1\textwidth]{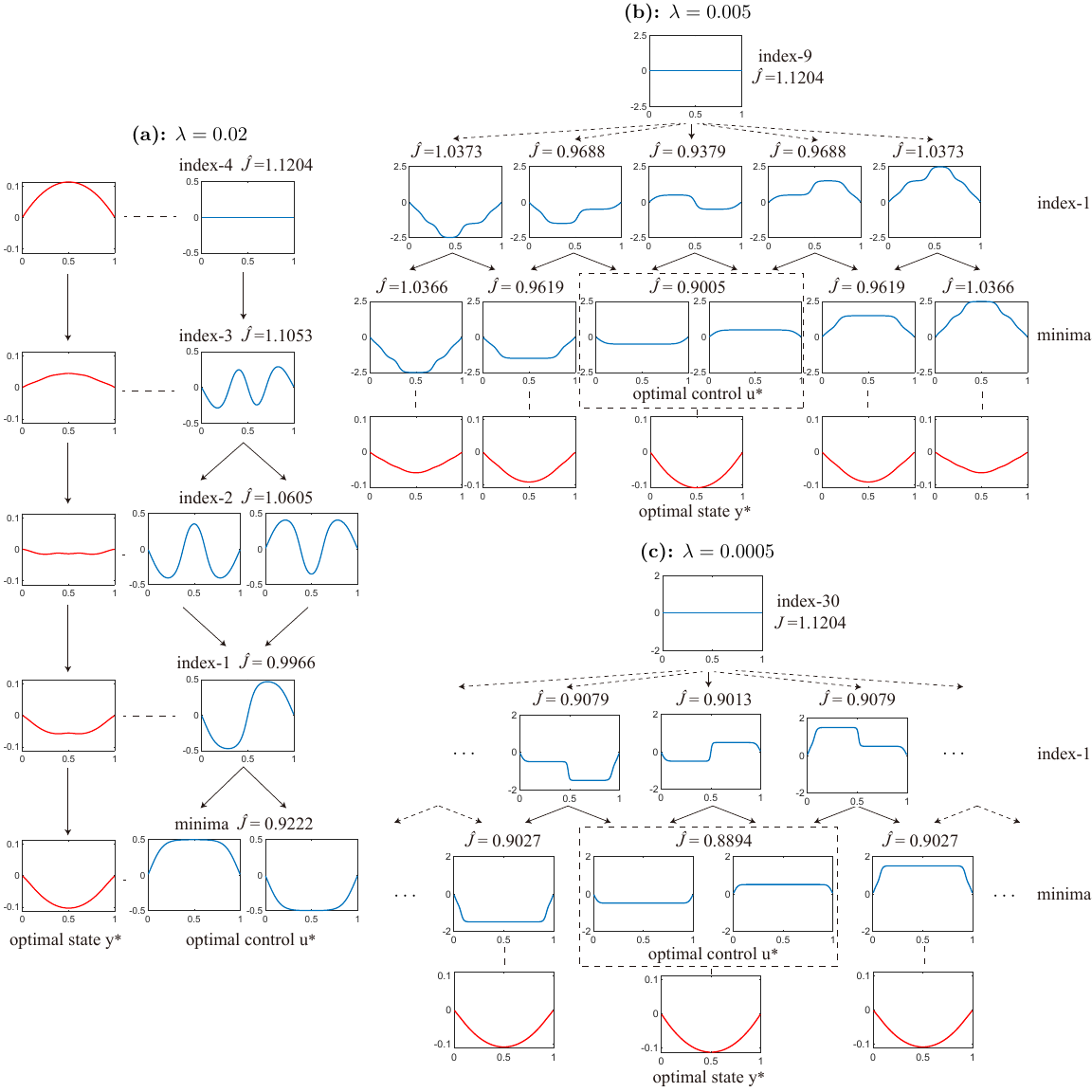}
    \caption{Control landscapes of one-dimensional problem with different $\lambda $. Each plot shows either a distribution of the control $u$ (blue curve) corresponding to a saddle point of the reduced cost functional $\hat{J}$, or the corresponding state $y$ (red curve). The optimal control and the corresponding optimal state are denoted as $u^*$ and $y^*$, respectively. }
     \label{fig:1D-SL}
\end{figure}

\vspace{2pt}
\underline{\textbf{Constructing control landscapes}}
For different $\lambda$, the control landscapes of (\ref{equ:nc-ocp})--(\ref{equ:se-nc-ocp})  are presented in Fig. \ref{fig:1D-SL}, and each control landscape starts from a trivial state $u\equiv 0$ to perform a downward search for lower-index saddle points.
When $\lambda = 0.02$, the control landscape of $\hat{J}$ and the  corresponding state landscape are shown,  where $u\equiv 0$ is an index-4 saddle point and two minima can be located along the transition pathway, which correspond to the lowest value of $\hat{J}$ such that they are most likely global minima (i.e. the optimal control). Furthermore, as the index of the control variable decreases, the corresponding state variable gets close to the target state $y_d$, as shown in the state landscape. These observations indicate that the control and state landscapes depict comprehensive patterns of control and state spaces of the optimal control, respectively, by means of stationary points, which provide rich information of optimal control and naturally locate its minima.

For the case $\lambda=0.005$, six minima are found, two of which correspond to the lowest value of $\hat{J}$ such that they are most likely optimal controls.
When $\lambda = 0.0005$, the stronger non-convexity results in a much richer control landscape with many minima. Nevertheless, the control landscape could still locate them. Note that the control landscape method does not depend on the initial guess of $u$ due to the employment of the information of high-index saddle points, which demonstrates its effectiveness and superiority. 

Another observation from Fig. \ref{fig:1D-SL} is that a decrement of $\lambda$ leads to a decrement of the lowest value of $\hat{J}$. A possible reason is that a decrement of $\lambda$ implies that the term $ \|u\|_{H^1(\Omega)}^2 $ plays a less dominant role in $\hat{J}$ such that the change of optimal control $u$ is admitted to be more rapid, cf.  Fig. \ref{fig:1D-SL}, leading to a better approximation $y$ to $y_d$ and thus a smaller value of $\hat{J}$.

\vspace{4pt}

\underline{\textbf{Comparison with Gradient and Newton's methods}}
For comparison with the control landscape method, we test the performance of the BB gradient method \cite{barzilaiTwoPointStepSize1988} and Newton's method in solving optimal control (\ref{equ:nc-ocp})--(\ref{equ:se-nc-ocp}). 

\begin{figure}[h] 
    \centering
    \includegraphics[width=1\textwidth]{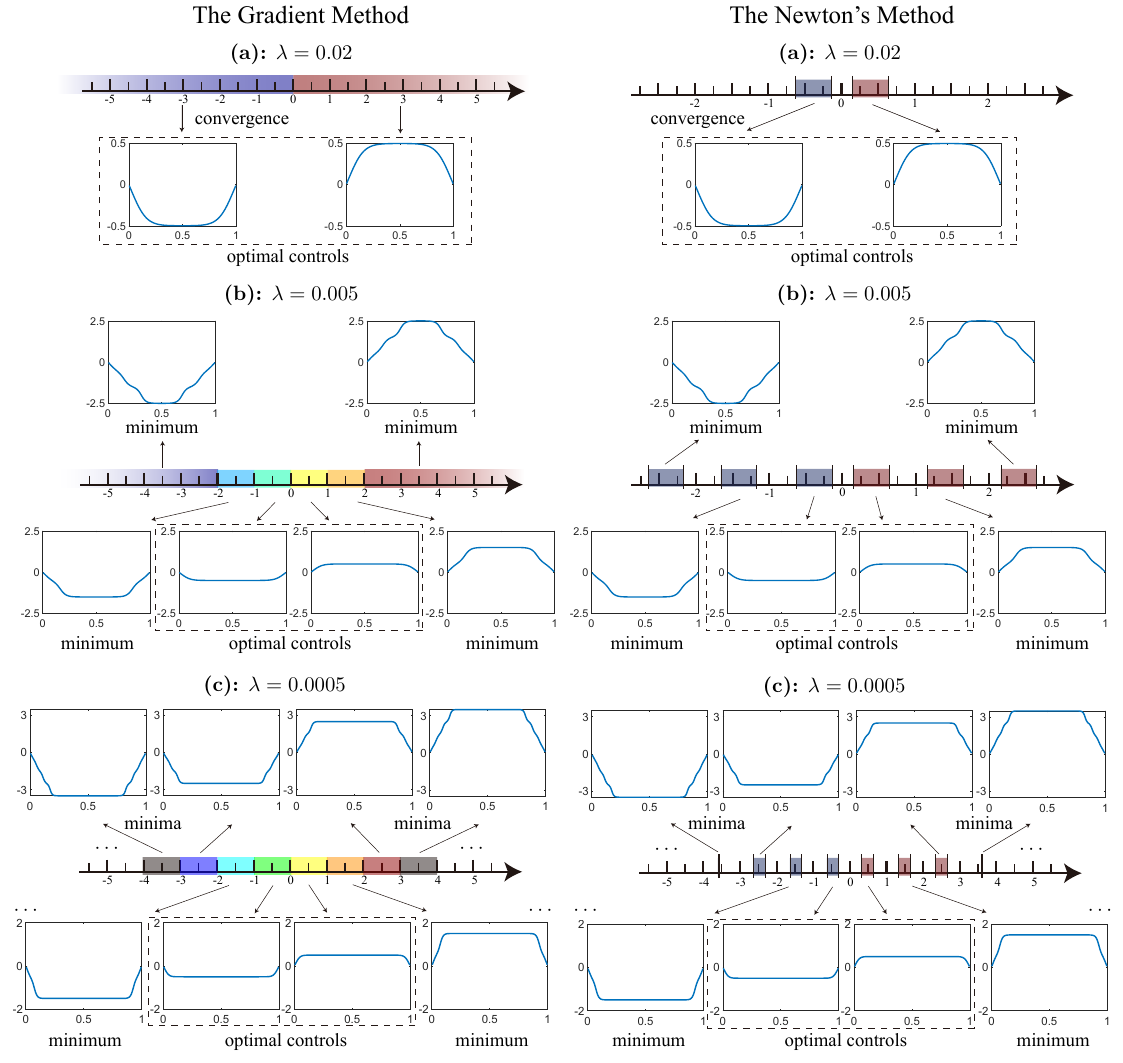}
    \caption{Domains of convergence of the gradient method and the Newton's method. Each colored segment represents an interval in which the constant initial value of $u$ could be selected such that the algorithm converges to the corresponding minima.}
     \label{fig:1D-NM}
\end{figure}

In both iterative methods, different constant initial values of the control variable $u$ are attempted over $[-10,10]$, that is, we test the initial value of $u$ as $-10+0.01j$ for $j=0,1,\ldots,2000$ in order, and the domains of convergence and the corresponding minima are presented in Fig. \ref{fig:1D-NM}. 
We observe that the convergence of both methods depends on the initial guesses of $u$, and in particular the algorithms could converge to optimal controls when the initial guesses are chosen from some narrow domains. As $\lambda$ decreases, the domains of convergence also shrink, which increases the difficulty of locating minima. These observations indicate that the effectiveness of locating minima by both methods rely on appropriate initial guesses. 
 
\subsection{Control landscapes for two-dimensional problems}
We consider examples on $\Omega=[0,1]^2$ with $a_{11}(x)=c(x)\equiv 1$, $ d(x,y) = 0$ and $y_d = 2$. A uniform partition is adopted with $h = 2^{-5}$, and different cases are considered:

\vspace{0.02in}

\begin{itemize}
    \item\textbf{Case 1} (symmetric, unconstrained): $g(x,u) = e^{-(u+1)^2} + e^{-(u-1)^2} - e^{-u^2}$; 
    \item\textbf{Case 2} (asymmetric, unconstrained): $ g(x,u) = e^{-(u+1)^2} + e^{-(u-1)^2} - (1+0.8u^3)e^{-u^2}$;
    \item\textbf{Case 3} (asymmetric, constrained): $g(x,u) = e^{-(u+1)^2} + e^{-(u-1)^2} - (1+0.8u^3)e^{-u^2}$
                with integral constraint  $\int_{\Omega} u dx = 0$.  
    \end{itemize}

\begin{figure}[h]
\centering
 \includegraphics[width=1\textwidth]{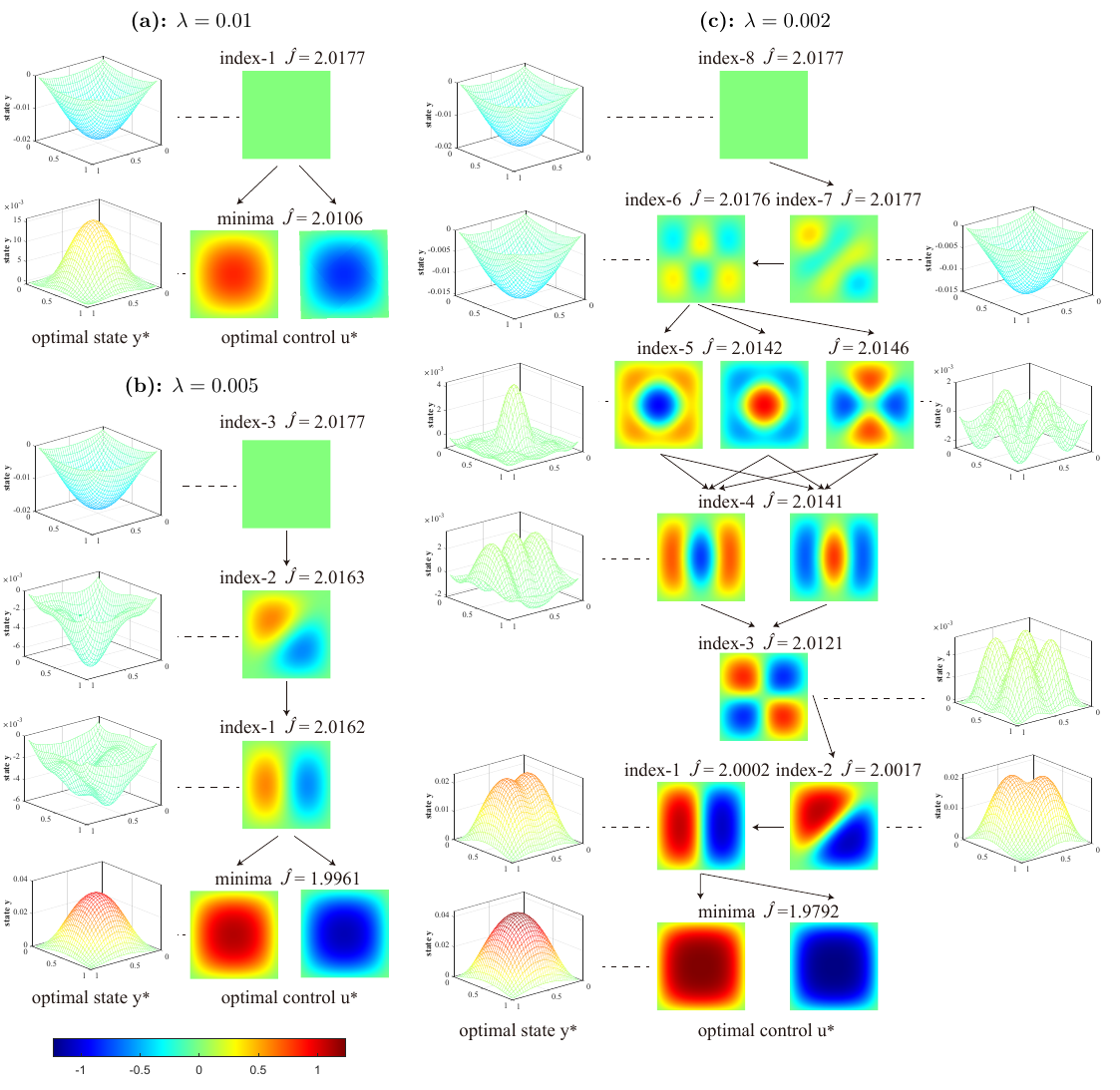}
\caption{Control landscapes of symmetric and unconstrained case 1 with $\lambda = 0.01,0.005, 0.002$. Each colored plot shows the distribution of the control $u$, and each three-dimensional surface depicts the corresponding state $y$. In this and all subsequent cases, saddle points that can be obtained from one another by rotations and sign reversals are regarded as identical.}
\label{fig:2DE1-SL}
\end{figure}

Here,  the additional zero mean constraint $\int_\Omega u\,dx=0$ is Case 3 is handled by introducing the projection operator $z\rightarrow z-\frac{1}{|\Omega|}\int_\Omega z dx$ in FHiSD, cf. the similar treatment for traditional HiSD in \cite{xuSolutionLandscapesDiblock2021}.
In Fig. \ref{fig:2DE1-SL}, 
we present control and state landscapes of $\hat{J}(u)$ for Case 1.
Similar to the one-dimensional case, a decrement of $\lambda$ leads to a more complex pattern, and during the evolution of the control variable along the transition pathways in control landscape, the state variable gradually evolves toward the desired state and eventually reaches a stable optimal configuration. These numerical results again indicate the effectiveness of the control landscape in depicting comprehensive patterns of the non-convex optimal control and locating minima. 

\begin{figure}[h]
    \centering
    \includegraphics[width=1\textwidth]{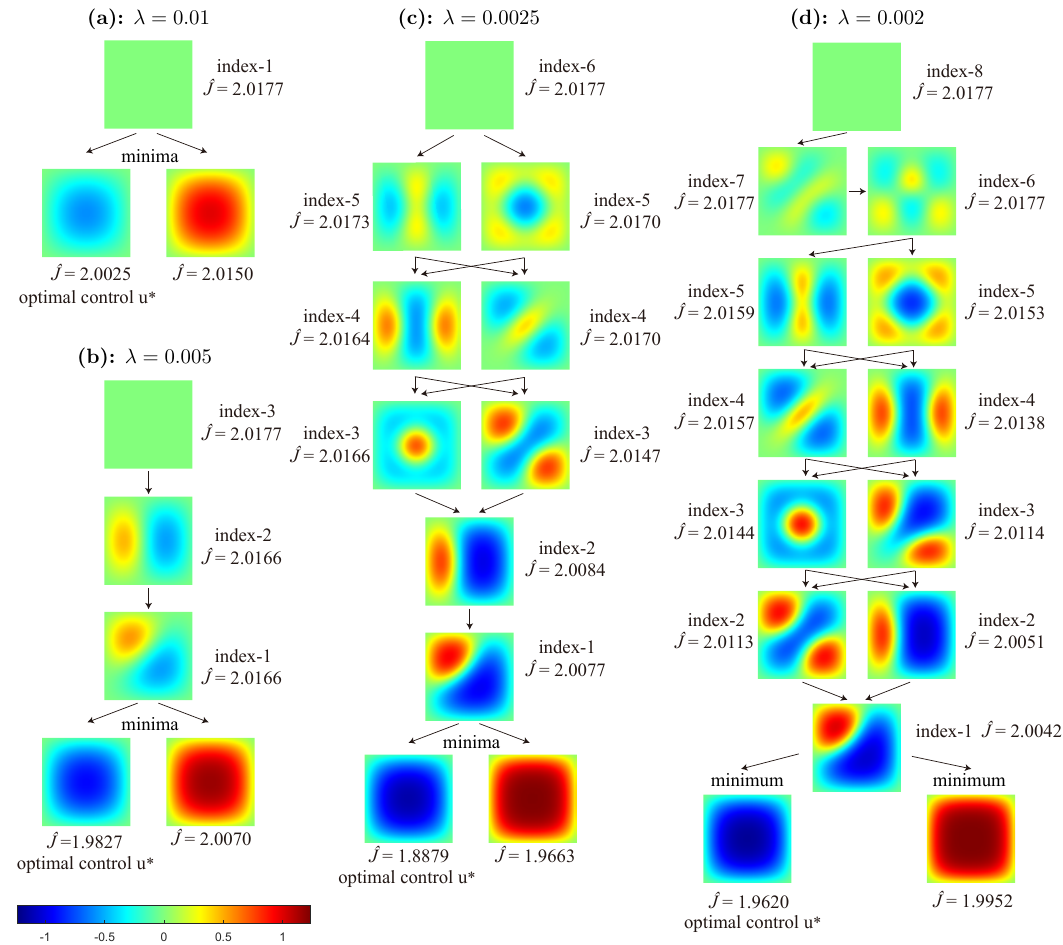}
    \caption{Control landscapes of asymmetric and unconstrained case 2 with $\lambda = 0.01,0.005,0.0025$, $0.002$.}
        \label{fig:2DE2-SL}
\end{figure}

In Fig. \ref{fig:2DE2-SL}, we consider the Case 2. Due to the broken symmetry, many deformations arise compared to the results in Fig. \ref{fig:2DE1-SL}. For instance, two local optimal controls are obtained for each case of $\lambda$, but only one of them attains the lowest value of $\hat{J}$. 
Another interesting and unintuitive phenomenon is that some high-index saddle points correspond to smaller values of $\hat{J}$ than those of the subsequent lower-index saddle points.
For instance, in Fig. \ref{fig:2DE2-SL} (d), an index-4 saddle point corresponds to  $\hat J=2.0138$, while an index-3 saddle point corresponds to $\hat{J}=2.0144$.
We thus reach the following conclusions or comments:
\begin{itemize}
\item
     A lower-index saddle point may not necessarily be more optimal than a higher-index one.
Moreover, asymmetric structures could further increase the non-convexity of the model and then break the energy-decreasing behavior along a transition pathway. Such phenomena are rarely encountered in other multiple solution problems and imply the necessity of computing multiple solutions of non-convex optimal control in order to find the optimal control.
\item In practical applications, the feasibility of implementing the optimal control should be considered.
If the optimal control is difficult to be realized in practice, one could still select other stationary points from the control landscape to balance the feasibility of the control and its optimality determined by $\hat{J}$.
\end{itemize}

\begin{figure}[h]
    \centering
    \includegraphics[width=1\textwidth]{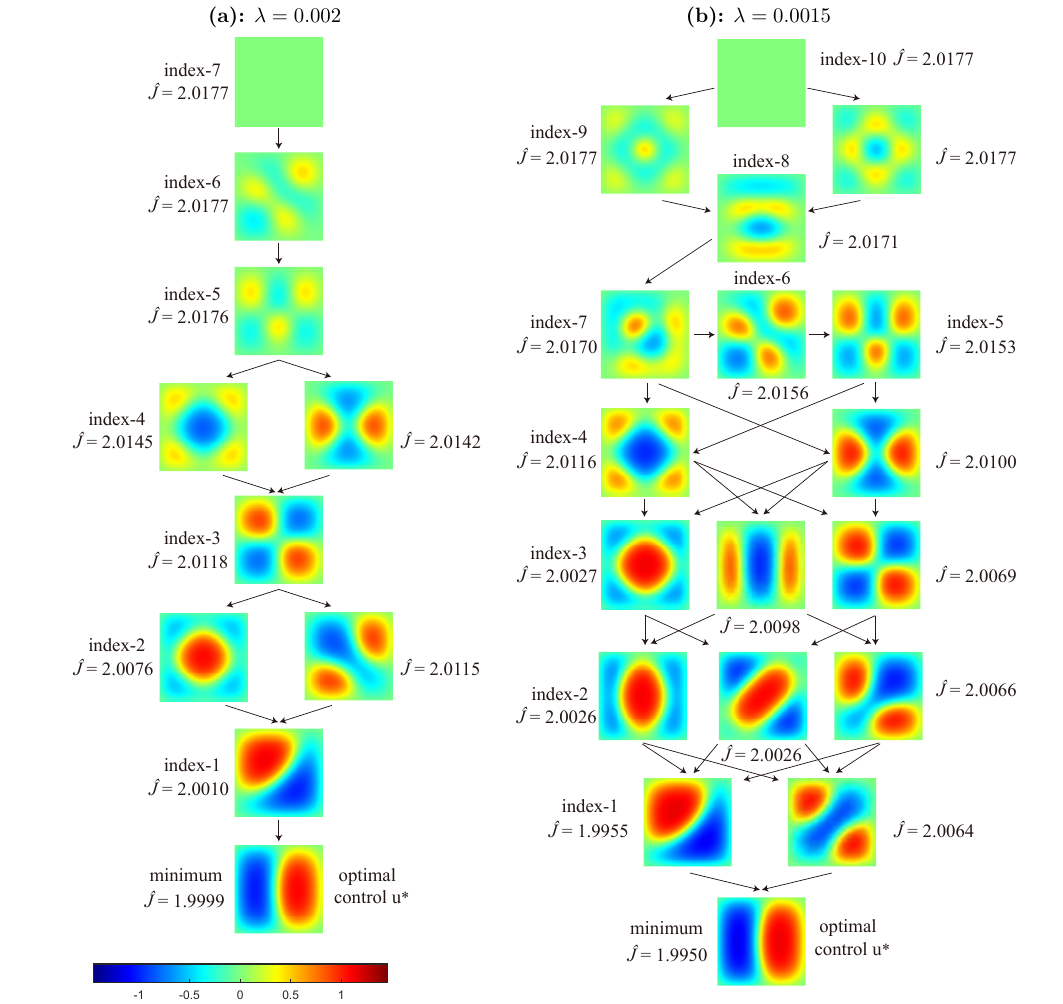}
    \caption{Control landscapes of asymmetric and constrained case 3 with $\lambda = 0.002, 0.0015$.}
        \label{fig:2DE3-SL}
\end{figure}

Finally, we consider the Case 3 in Fig. \ref{fig:2DE3-SL}, which shows that under the integral constraint the optimal control  problem  admits a unique optimal control.  The distribution of optimal control has a sharp interface due  to the integral constraint, which is different from the unconstrained case. Numerical results again demonstrate the effectiveness of the control landscape in solving constrained optimal control and the applicability of the idea of imposing the constraint on FHiSD. 
 
\section{Concluding remarks}

In this work, we solve the non-convex optimal control by using FHiSD to construct the control landscape, which could locate multiple stationary points with less requirements on initial guesses. Numerical results indicate the strength of the control landscape method in solving the non-convex optimal control. 
It is worth mentioning that though the current work focuses on an optimal control of elliptic equation, the proposed idea has a great extensibility and could be widely extended to other optimal control and optimization problems. 
For instance, optimal control with different cost functionals subject to different constrained equations and admissible sets could be considered, though some details still need further exploration, e.g. how to ensure that the primary variable of FHiSD evolves in the admissible set is not a trivial issue. 
Furthermore, the control landscape contains much more information than the location of minima, and numerical results reveal rich solution structures of the non-convex optimal control. Whether such information has physical meanings and potential applications remains to be further studied.

\backmatter

\bmhead{Acknowledgements}

This work was partially supported by the National Key R\&D Program of China (Grant Nos. 2023YFA1009200 and 2023YFA1008903), 
the National Natural Science Foundation of China (Grant Nos. 12225102, T2321001, 12288101, 12301555 and 12131014), 
the Natural Science Foundation of Shandong Province (Grant No. ZR2025QB01) and   
the Taishan Scholars Program of Shandong Province (Grant No. tsqn202306083).

\section*{Appendix: Derivations of \eqref{equ:df}--\eqref{equ:dddf}}
\label{app:derivative-proof} 
Based on Theorem \ref{thm:diff}, the reduced cost functional $\hat{J}(u)$ is three times Fr\'echet differentiable on $H^1_0(\Omega)$.
Now we intend to characterize the derivatives of $\hat J(u)$.

For $u, z_1\in H_0^1(\Omega)$, we invoke (\ref{equ:ae-nc-ocp}) to obtain
\begin{equation*}
\begin{aligned}
\langle \hat{J}'(u),z_1\rangle
&= \frac{d}{dt}  \hat{J}(u+t z_1) |_{t=0} = \int_{\Omega} \{(y-y_d)G'(u)z_1 + \lambda (u z_1 + \nabla u \nabla z_1 )\} dx\\
&= \int_{\Omega} (g_u(x,u) p+ \lambda u) z_1 dx + \int_{\Omega}  \lambda \nabla u \nabla z_1 dx.
\end{aligned}
\end{equation*} 

\vspace{-0.1in}

Similarly, for $u,z_1, z_2 \in H_0^1(\Omega)$, using \eqref{equ:ae-nc-ocp}, \eqref{equ:linearized-state} and \eqref{equ:second-linearized-state}, we have
\begin{equation*}
\begin{aligned}
\langle \hat{J}''(u)z_1, z_2\rangle
&= \frac{d}{dt} [ \langle \hat{J}'(u+t z_2),z_1\rangle ] |_{t=0} \\
&\hspace{-0.7in}= \int_{\Omega} [ (y-y_d)\delta y_{z_1 z_2} + \delta y_{z_1} \delta y_{z_2}+ \lambda (z_1 z_2 + \nabla z_1 \nabla z_2 ) ]dx\\
        & \hspace{-0.7in}= \int_{\Omega} [ (\mathcal{A} p + d_y(x,y)p)\delta y_{z_1 z_2} + \delta y_{z_1}\delta y_{z_2} ] dx + \lambda\!\int_{\Omega} (z_1 z_2 + \nabla z_1 \nabla z_2 ) dx \\
        &\hspace{-0.7in} = \int_{\Omega}( \lambda z_1 + p\, g_{uu}(x,u) z_1) z_2 \,dx + \lambda \int_{\Omega} \nabla z_1 \nabla z_2 \,dx  +\int_{\Omega} (1 - p\, d_{yy}(x,y))  \delta y_{z_1} \delta y_{z_2} \,dx.
\end{aligned}
\end{equation*} 

For the third derivative, let $z_3\in H_0^1(\Omega)$, $ \delta y_{z_3}=G'(u)[z_3]$, and let  $\delta p_3$ denote the derivative of the adjoint state in the direction $z_3$, which is determined by
\begin{equation}\label{equ:2d-adjoint}
    \begin{cases}
    \mathcal A\delta p_3+d_y(x,y)\delta p_3
    =
    \bigl(1-p\,d_{yy}(x,y)\bigr)\delta y_{z_3}
    & \text{in }\Omega,\\
    \delta p_3=0
    & \text{on }\Gamma.
    \end{cases}
\end{equation}
For simplicity, set $q:=1-p\,d_{yy}(x,y)$.
Then, differentiating $\langle \hat J''(u) z_1, z_2\rangle$ in the direction $z_3$ yields
\begin{equation*}
\begin{aligned}
\langle \hat J'''(u) z_1 z_2, z_3 \rangle 
&=
\frac{d}{dt}
\left[
\langle \hat J''(u+t z_3) z_1, z_2\rangle
\right]_{t=0} \\
&=
\bigl(\delta p_3\,g_{uu}(x,u)z_1 z_2,1\bigr)
+\bigl(p\,g_{uuu}(x,u) z_1 z_2 z_3,1\bigr)\\
&\quad
+\bigl(q\,\delta y_{z_1 z_3},\delta y_{z_2}\bigr)
+\bigl(q\,\delta y_{z_1},\delta y_{z_2 z_3}\bigr)\\
&\quad
-\bigl(\delta p_3\,d_{yy}(x,y)\delta y_{z_1}\delta y_{z_2},1\bigr)
-\bigl(p\,d_{yyy}(x,y)\delta y_{z_1}\delta y_{z_2}\delta y_{z_3},1\bigr).
\end{aligned}
\end{equation*}
Using \eqref{equ:second-linearized-state} and \eqref{equ:2d-adjoint}, we have
\begin{equation*}
\bigl(\delta p_3, \,g_{uu}(x,u)z_1 z_2 - d_{yy}(x,y)\delta y_{z_1}\delta y_{z_2}\bigr) 
\!=\!
\bigl(\delta p_3,(\mathcal A+d_y(x,y))\delta y_{z_1 z_2}\bigr)
\!=\! 
\bigl(q\,\delta y_{z_1 z_2},\delta y_{z_3}\bigr),
\end{equation*}
which leads to
\begin{equation*}
    \begin{aligned}
        \langle \hat J'''(u)z_1 z_2, z_3 \rangle 
        &=
        \bigl(q\delta y_{z_1 z_2},\delta y_{z_3}\bigr) 
        +
        \bigl(q\delta y_{z_1 z_3},\delta y_{z_2}\bigr) 
        + 
        \bigl(q\delta y_{z_2 z_3},\delta y_{z_1}\bigr)\\
        &\quad + 
        \bigl(p\,g_{uuu}(x,u) z_1 z_2 z_3,1\bigr)
        -
        \bigl(p\,d_{yyy}(x,y)\delta y_{z_1}\delta y_{z_2}\delta y_{z_3},1\bigr).
    \end{aligned}
\end{equation*}
Thus we complete the derivation of \eqref{equ:df}, \eqref{equ:ddf} and \eqref{equ:dddf}.

\nocite{scmrefs}
\bibliography{references}

\end{document}